\documentclass[letter,10pt]{amsart}
\usepackage[english]{babel} 
\usepackage{multicol}
\usepackage{graphicx}
\usepackage{amssymb}
\usepackage[all]{xy}
\usepackage{amsmath,amsfonts}
\usepackage{enumerate}
\usepackage{array}
\usepackage{tikz}
\input xy
\xyoption{all}

\usepackage{xcolor}

\newtheorem{teo}{Theorem}[section]
\newtheorem{cor}[teo]{Corollary}
\newtheorem{lema}[teo]{Lemma}
\newtheorem{prop}[teo]{Proposition}
\newtheorem{defi}[teo]{Definition}

\newtheorem{Rec}[teo]{Recall}
\newtheorem{Rem}[teo]{Remark}

\numberwithin{equation}{section}

\title{Third homology of $\mathrm{SL}_{2}$ over Number fields: The norm-Euclidean quadratic imaginary case.}
\author{Rodrigo Cuitun Coronado}
\address{School of Mathematics and Statistics, University College Dublin}
\email{rodrigo.cuituncoronado@ucdconnect.ie}
\date{\today}

\begin{document}

\begin{abstract}
In the article \emph{The third homology of} $SL_{2}(\mathbb{Q})$ (\cite{ArticlehomologyQprint}), Hutchinson determined the structure of
$H_{3}\left(\mathrm{SL}_{2}(\mathbb{Q}),\mathbb{Z}\left[\frac{1}{2}\right]\right)$
by expressing it in terms of $K_{3}^{\mathrm{ind}}(\mathbb{Q})\cong
\mathbb{Z}/24$ and the \emph{scissor congruence group} of the residue field
$\mathbb{F}_{p}$ with $p$ a prime number. In this paper, we develop further
the properties of \emph{the refined scissors congruence group} in order to
extend this result to the case of imaginary quadratic number fields whose ring of
integers is a Euclidean domain with respect to the norm.

\end{abstract}
\keywords{Group Homology, Special linear group, $K$-theory, Number Fields}

\maketitle

\section{Introduction}

In the recent paper \emph{The third homology of} $SL_{2}(\mathbb{Q})$ (\cite{ArticlehomologyQprint}), Hutchinson calculates the structure of $H_{3}(\mathrm{SL}_{2}(\mathbb{Q}),\mathbb{Z})$ by calculating the surjective homomorphism $H_{3}(\mathrm{SL}_{2}(\mathbb{Q}),\mathbb{Z})\rightarrow K_{3}^{\mathrm{ind}}(\mathbb{Q})$. 
He asked whether the structure theorem proved in his paper generalize to number fields and even global fields. We answer this question in this paper, the case of imaginary quadratic number fields whose ring of integers is a Euclidean domain with respect to the norm of the field (see section \ref{Section_proof_of_the_main_theorem} below).

\medskip

Let $F$ be a field, there is a natural induced map from $H_{3}(\mathrm{SL}_{2}(F), \mathbb{Z})$ to the indecomposable quotient,$K_{3}^{\mathrm{ind}}(F)$, of the third Quillen $K$-group of $F$. It can be shown this map is surjective (see \cite{HutchinsonLiqunTao}). In \cite{sah:dupont}, Johan Dupont and Chih-Han Sah showed that this map is an isomorphism for the case of $F=\mathbb{C}$ (and more generally $F^{\times}=(F^{\times})^{2}$). Furthermore, in \cite{ChinHanSah}, Chih-Han Sah proved this result for the case of $\mathbb{R}$. However, for general fields, the map $H_{3}(\mathrm{SL}_{2}(F), \mathbb{Z})\rightarrow K_{3}^{\mathrm{ind}}(F)$ has a large kernel (\cite{ArticleBloch}). In fact, $H_{3}(\mathrm{SL}_{2}(F),\mathbb{Z})$ is naturally a module over the integral group ring $R_{F}:=\mathbb{Z}[F^{\times}/(F^{\times})^2]$ of the group of square classes of the field $F$ and the above homomorphism is a homomorphism of $R_{F}$-modules if $K_{3}^{\mathrm{ind}}(F)$ is given the trivial module structure. In (\cite{ArticleBloch}), Hutchinson shows that the action of $R_{F}$ on $H_{3}(\mathrm{SL}_{2}(F), \mathbb{Z})$ is highly nontrivial for fields with (many) discrete valuations. When $F^{\times}$ acts nontrvially on $H_{3}(\mathrm{SL}_{2}(F),\mathbb{Z})$, the map  $H_{3}(\mathrm{SL}_{2}(F),\mathbb{Z})\rightarrow K_{3}^{\mathrm{ind}}(F)$ has a nontrivial kernel. Mirzaii has shown (\cite{mirzaii:third}) that the kernel of the induced homomorphism $H_{3}(\mathrm{SL}_{2}(F),\mathbb{Z})_{F^{\times}}\rightarrow K_{3}^{\mathrm{ind}}(F)$  consists at most of $2$-primary torsion. In other words $\mathrm{Ker}\left(H_{3}\left(\mathrm{SL}_{2}(F), \mathbb{Z}\left[\frac{1}{2}\right]\right)\rightarrow K_{3}^{\mathrm{ind}}(F)\left[\frac{1}{2}\right]\right)$ is $\mathcal{I}_{F}H_{3}(\mathrm{SL}_{2}(F), \mathbb{Z})$ (up to some $2$-torsion), where $\mathcal{I}_{F}$ is the augmentation ideal of $R_{F}$.

\medskip

For $F=\mathbb{Q}(\sqrt{-m})$  where $m\in \{1,2,3,7,11\}$, the main proposition of this article (see Proposition \ref{Main_Proposition_Isomorphism}) describes the structure of $H_{0}\left(\mathcal{O}^{\times}_{F}, \mathcal{I}_{F}H_{3}\left(\mathrm{SL}_{2}(F), \mathbb{Z}\left[\frac{1}{2}\right]\right)\right)$ as a $R_{\mathbb{Q}(\sqrt{-m})}$-module. The proposition states that this module - via a natural residue homomorphism $S_{p}$ - is isomorphic to the direct sum of the \emph{scissor congruence group} of the fields $k(\upsilon_{p})$, $\mathcal{P}(k(\upsilon_{p}))\left[\frac{1}{2}\right]$, where the sum is over all the primes of the field $F$(see Section \ref{Section_Fields_with_a_valuation} and Section \ref{Section_proof_of_the_main_theorem} for more details). The square class of $-1$ acts trivially on each factor, while the square class of $p$ acts by multiplication by $-1$ on the corresponding factor $\mathcal{P}(k(\upsilon_{p}))$. When $m\not=1$, it follows that as an abelian group

\begin{center}
$H_{3}\left(\mathrm{SL}_{2}(\mathbb{Q}(\sqrt{-m})),\mathbb{Z}\left[\frac{1}{2}\right]\right)\cong K_{3}^{\mathrm{ind}}(\mathbb{Q}(\sqrt{-m}))\left[\frac{1}{2}\right]\oplus\left(\displaystyle\bigoplus_{p\in\Pi}\mathcal{P}(k(\upsilon_{p}))\left[\frac{1}{2}\right]\{p\}\right)$
\end{center}

where $\Pi$ denotes the set of primes (since $\mathcal{O}_{F}^{\times}$ acts trivially in this case) and the structure of indecomposable $K_{3}$ over number fields and the \emph{scissor congruence group} of finite fields are known (see \cite{KtheoryCharlesWeibel} and \cite{ArticleBlochWignerComplex} respectively).



\medskip

One of the main tool we use is the description of $H_{3}\left(\mathrm{SL}_{2}(F),\mathbb{Z}\left[\frac{1}{2}\right]\right)$ in terms of the \emph{refined scissors congruence group}. The \emph{scissors congruence group} $\mathcal{P}(F)$ of a field $F$ was introduced by Dupont and Sah in their study of the Hilbert's third problem in hyperbolic $3$-space (\cite{sah:dupont}). It is an abelian group defined by an explicit presentation (see section \ref{Section_Bloch_group_of_fields}) and it was shown for the case $F=\mathbb{C}$ (and more generally when $F$ is algebraically closed) by the authors that $K_{3}^{\mathrm{ind}}(\mathbb{C})\cong H_{3}(\mathrm{SL}_{2}(\mathbb{C}),\mathbb{Z})$. Soon afterwards, Suslin (\cite{Suslin1}) proved that, for any infinite field $F$ it is closely connected with $K_{3}^{\mathrm{ind}}(F)$, the indecomposable $K_{3}$ of $F$ (see Theorem \ref{Suslinmaintheorem}). As noted above, $H_{3}(\mathrm{SL}_{2}(F),\mathbb{Z})$ is naturally a module over the integral group ring $R_{F}:=\mathbb{Z}[F^{\times}/(F^{\times})^2]$ of the group of square classes of the field $F$. The \emph{refined scissors congruence group} of the field $F$ - introduced by Hutchinson in \cite{ArticleBlochWignerComplex}- is defined by generators and relations analogously to the \emph{scissors congruence group} but as  a module over $R_{F}$ (see section \ref{Section_Refined_Bloch_group_of_fields}). In \cite{ArticleBlochWignerComplex}, it is elucidated the relation of the \emph{refined scissors congruence group} $\mathcal{RP}(F)$ to $H_{3}(\mathrm{SL}_{2}(F),\mathbb{Z})$ and $K_{3}^{\mathrm{ind}}(F)$ (see Theorem \ref{TheoremBlochWignerComplex} for a precise statement). Our starting point in this article is essentially an isomorphism $\mathcal{I}_{F}H_{3}\left(\mathrm{SL}_{2}(F)), \mathbb{Z}\left[\frac{1}{2}\right]\right)\cong \mathcal{I}_{F}\mathcal{RP}_{+}(F)\left[\frac{1}{2}\right]$ where $\mathcal{RP}_{+}(F)$ is a certain quotient of $\mathcal{RP}(F)$ (see Proposition \ref{PropsitionRP+(F)andH3(SL2)}).

\subsection{Layout of this article}

In Section \ref{Section_Bloch_group_of_fields}, we review some known results about the scissor congruence groups and their relation to the third homology of $\mathrm{SL}_{2}$ of fields.

In Section \ref{Section_Characters}, we review the character-theoretic local-global principle for modules over a group ring (\cite{Articlediscretevaluation}). We recall some applications to \emph{Scissor congruence group} (\cite{ArticlehomologyQprint}), and generalize some of these results.

In Section \ref{Section_Fields_with_a_valuation}, we review some algebraic properties of the refined scissor congruence groups over fields $F$ which have a discrete valuation.

In Section \ref{Section_Quadratic_fields}, we review relevant facts about quadratic number fields.

Section \ref{Section_proof_of_the_main_theorem} contains the proof of the main proposition (Proposition \ref{Main_Proposition_Isomorphism}) for the field $\mathbb{Q}(\sqrt{-m})$, where $m\in \{1,2,3,7,11\}$. 


\subsection{Some notation}

For a Field $F$, we let $F^{\times}$ denote the group of units of $F$. For $x\in F^{\times}$ we will let $\langle x\rangle\in F^{\times}/(F^{\times})^{2}$ denote the corresponding square class. Let $R_{F}$ denote the integral group ring $\mathbb{Z}\left[F^{\times}/(F^{\times})^{2}\right]$ of the group $F^{\times}/(F^{\times})^{2}$. We will use the notation $\langle\langle x\rangle\rangle$ for the basis elements, $\langle x\rangle-1$ of the augmentation ideal $\mathcal{I}_{F}$ of $R_{F}$. For any $a\in F^{\times}$, we will let $p^{a}_{+}$ and $p^{a}_{-}$ denote the elements $1+\langle a\rangle$ and $1-\langle a\rangle$ in $R_{F}$ respectively. For any abelian group $G$ we will let $G\left[\frac{1}{2}\right]$ denote $G\otimes\mathbb{Z}\left[\frac{1}{2}\right]$.\newline

\subsection{$H_{\bullet}(\mathrm{SL}_{2}(F), \mathbb{Z})$ is a $R_{F}$-module}

Let us recall the group extension

\begin{center}
$\xymatrix{1\ar[r]& \mathrm{SL}_{n}(F)\ar[r]& \mathrm{GL}_{n}(F)\ar[r]&F^{\times}\ar[r]&1}$
\end{center}

induces an action - by conjugation - of $F^{\times}$ on  the homology group $H_{\bullet}(\mathrm{SL}_{n}(F), \mathbb{Z})$. Since the determinant of a scalar matrix is an $n$-th power, the subgroup $(F^{\times})^{n}$ acts trivially. In particular, the groups
$H_{\bullet}(\mathrm{SL}_{2}(F), \mathbb{Z})$ are modules over the integral group ring $R_{F}:=\mathbb{Z}[F^{\times}/(F^{\times})^2]$.

\section{Bloch Groups of fields.}\label{Section_Bloch_group_of_fields}

In this section we will recall the definition and applications of the classical pre-Bloch group $\mathcal{P}(F)$ and the refined pre-Bloch group $\mathcal{RP}(F)$.

\subsection{Classical Bloch Group $\mathcal{B}(F)$.}
For a field $F$, with at least $4$ elements, the \emph{pre-Bloch group} or \emph{Scissors congruence group}, $\mathcal{P}(F)$, is the group generated by the elements $[x]$, with $F^{\times}\setminus\{1\}$, subject to the relations

\begin{center}
$R_{x,y}:\; 0=[x]-[y]+[y/x]-[(1-x^{-1})/(1-y^{-1})]+[(1-x)/(1-y)]$, \;\; $x\not=y$.
\end{center}

Let $S^{2}_{\mathbb{Z}}(F^{\times})$ denote the group

\begin{center}
$\frac{F^{\times}\otimes_{\mathbb{Z}}F^{\times} }{\langle x\otimes y +y\otimes x\; |\;x,y\in F^{\times}\rangle}$
\end{center}

and denote by $x\circ y$ the image of $x\otimes y$ in $S^{2}_{\mathbb{Z}}(F^{\times})$. The map

\begin{center}
$\lambda: \mathcal{P}(F)\rightarrow S^{2}_{\mathbb{Z}}(F^{\times})$,  \;\; $[x]\mapsto (1-x)\circ x$
\end{center}


is well-defined, and the \emph{Bloch group} of $F$, $\mathcal{B}(F)\subset\mathcal{P}(F)$, is defined to be the kernel of $\lambda$.

\smallskip

For the fields $\mathbb{F}_{2}$ and $\mathbb{F}_{3}$ the following definitions allows us to include these fields in the statements of some of our results:

\smallskip

$\mathcal{P}(\mathbb{F}_{2})=\mathcal{B}(\mathbb{F}_{2})$ is a cyclic group of order $3$ with generator denoted $C_{\mathbb{F}_{2}}$.

\smallskip

$\mathcal{P}(\mathbb{F}_{3})$ is cyclic of order $4$ with generator $[-1]$. $\mathcal{B}(\mathbb{F}_{3})$ is the subgroup generated by $2[-1]$.

The Bloch group $\mathcal{B}(F)$ of a general field $F$ is of interest because of the following result of Suslin on $K_{3}^{\mathrm{ind}}$:

\begin{teo}\cite[Theorem 5.2]{Suslin1}\label{Suslinmaintheorem}
Let $F$ be an infinite field, then there is a short exact sequence

\begin{center}
$\xymatrix{0\ar[r]&\mathrm{Tor}_{1}^{\mathbb{Z}}(\widetilde{\mu_{F},\mu_{F}})\ar[r]&K_{3}^{\mathrm{ind}}(F)\ar[r]&\mathcal{B}(F)\ar[r]&0}$
\end{center}

where $\mathrm{Tor}_{1}^{\mathbb{Z}}(\widetilde{\mu_{F},\mu_{F}})$  is the unique nontrivial extension of $\mathrm{Tor}_{1}^{\mathbb{Z}}(\mu_{F},\mu_{F})$ by $\mathbb{Z}/2$ when $\mathrm{Char}(F)\not=2$ (and $\mathrm{Tor}_{1}^{\mathbb{Z}}(\widetilde{\mu_{F},\mu_{F}})=\mathrm{Tor}_{1}^{\mathbb{Z}}(\mu_{F},\mu_{F})$ if $\mathrm{Char}(F)=2$) .
\end{teo}

\begin{flushright}
    $\Box$
\end{flushright}

\subsection{The refined Bloch Group $\mathcal{RB}(F)$.}\label{Section_Refined_Bloch_group_of_fields}

The \emph{refined pre-Bloch group} $\mathcal{RP}(F)$, of a field $F$ which has at least $4$ elements, is the $R_{F}$-module with generators $[x]$, $x\in F^{\times}$ subject to the relations $[1]=0$ and

\begin{center}
$S_{x,y}:\; 0=[x]-[y]+\langle x\rangle[y/x]-\langle x^{-1}-1\rangle[(1-x^{-1})/(1-y^{-1})]+\langle1-x\rangle[(1-x)/(1-y)]$, \;\; $x,y\not=1$.
\end{center}

From the definitions, it follows that $\mathcal{P}(F)=(\mathcal{RP}(F))_{F^{\times}}$. Let $\Lambda=(\lambda_{1},\lambda_{2})$ be the $R_{F}$-module homomorphism

\begin{center}
$\mathcal{RP}(F)\rightarrow \mathcal{I}_{F}^{2}\oplus S^{2}_{\mathbb{Z}}(F^{\times})$
\end{center}

where $S^{2}_{\mathbb{Z}}(F^{\times})$ has the trivial $R_{F}$-module structure, $\lambda_{1}:\mathcal{RP}(F)\rightarrow \mathcal{I}^{2}_{F}$ is the map $[x]\mapsto \langle\langle1-x\rangle\rangle\langle\langle x\rangle\rangle$ and $\lambda_{2}$ is the composite

\begin{center}
$\xymatrix{ \mathcal{RP}(F)\ar[r]&\mathcal{P}(F)\ar[r]^{\lambda} &S^{2}_{\mathbb{Z}}(F^{\times})}$,
\end{center}


The \emph{refined Bloch group} of $F$ is the module

\begin{center}
$\mathcal{RB}(F):=\mathrm{Ker}\left(\Lambda: \mathcal{RP}(F)\rightarrow \mathcal{I}_{F}^{2}\oplus S^{2}_{\mathbb{Z}}(F^{\times}  )\right)$.
\end{center}

Furthermore, the \emph{refined scissor congruence group} of $F$ is the $R_{F}$-module

\begin{center}
$\mathcal{RP}_{1}(F):=\mathrm{Ker}(\lambda_{1}:\mathcal{RP}(F)\rightarrow\mathcal{I}^{2}_{F})$.
\end{center}

Thus $\mathcal{RB}(F)=\mathrm{Ker}(\lambda_{2}:\mathcal{RP}_{1}(F)\rightarrow S^{2}_{\mathbb{Z}}(F^{\times}))$.

\bigskip

The refined Bloch group is of interest because of the following result on the third homology of $\mathrm{SL}_{2}$ over a general field $F$:

\begin{teo}\cite[Theorem 4.3]{ArticleBlochWignerComplex}\label{TheoremBlochWignerComplex}
Let $F$ be a field with at least $4$ elements.

\begin{itemize}
\item[(1)]If $F$ is infinite, there is a natural complex of $R_{F}$-modules

\begin{center}
$\xymatrix{0\ar[r]&\mathrm{Tor}_{1}^{\mathbb{Z}}(\mu_{F},\mu_{F})\ar[r]&H_{3}(\mathrm{SL}_{2}(F),\mathbb{Z})\ar[r]& \mathcal{RB}(F)\ar[r]&0}$,
\end{center}

which is exact except at the middle term where the homology is annihilated by $4$.

\item[(2)]If $F$ is finite of odd characteristic, there is a complex

\begin{center}
$\xymatrix{0\ar[r]&H_{3}(B,\mathbb{Z})\ar[r]&H_{3}(\mathrm{SL}_{2}(F),\mathbb{Z})\ar[r]& \mathcal{RB}(F)\ar[r]&0}$,
\end{center}

which is exact except at the middle term, where the homology has order $2$.

\item[(3)]If $F$ is finite of characteristic $2$, there is a exact sequence

\begin{center}
$\xymatrix{0\ar[r]&H_{3}(B,\mathbb{Z})\ar[r]&H_{3}(\mathrm{SL}_{2}(F),\mathbb{Z})\ar[r]& \mathcal{RB}(F)\ar[r]&0}$,
\end{center}
\end{itemize}

\end{teo}

\begin{flushright}
    $\Box$
\end{flushright}

Now for $x\in F^{\times}$, we define the following elements of $\mathcal{RP}(F)$

\begin{center}
$\psi_{1}(x):=[x]+\langle-1\rangle[x^{-1}]$ \;\; and\;\; $\psi_{2}(x):=\left\{
                                                                       \begin{array}{ll}
                                                                         \langle x^{-1}-1\rangle[x]+\langle1-x\rangle[x^{-1}], & \hbox{$x\not=1$;} \\
                                                                         0, & \hbox{$x=1$.}
                                                                       \end{array}
                                                                     \right.
$.
\end{center}

From the definitions of the the elements $\psi_{i}(x)$, we get that $\langle-1\rangle\psi_{i}(-1)=\psi_{i}(-1)$ for $i\in\{1,2\}$.
We define $\widetilde{\mathcal{RP}}(F)$ to be $\mathcal{RP}(F)$ modulo the submodule generated by the elements $\psi_{1}(x)$, $x\in F^{\times}$.

\smallskip

In the section $3.2$ of \cite{Articlediscretevaluation}, it is shown that the elements

\begin{center}
$C(x)=[x]+\langle-1\rangle[1-x]+ \langle\langle1-x\rangle\rangle\psi_{1}(x) \in \mathcal{RP}(F)$
\end{center}

are constant for a field with at least 4 elements; i.e. $C(x)=C(y)$ for all $x,y\in F^{\times}$. Therefore we have the following definition.

\begin{defi}
Let $F$ be a field with at least $4$ elements. We will denote by $C_{F}$ the common value of the expression $C(x)$ for $x\in F\setminus\{0,1\}$; i.e.

\begin{center}
$C_{F}:=[x]+\langle-1\rangle[1-x]+ \langle\langle1-x\rangle\rangle\psi_{1}(x)$ \;\;in\;\; $\mathcal{RP}(F)$.
\end{center}

\end{defi}

For the field with $2$ and $3$ elements the following definitions allow us to include these fields in the statements of some results:

\smallskip

$\mathcal{RP}(\mathbb{F}_{2})=\mathcal{RP}_{1}(\mathbb{F}_{2})=\mathcal{RB}(\mathbb{F}_{2})=\mathcal{P}(\mathbb{F}_{2})$; i.e. it a cyclic module generated by $C_{\mathbb{F}_{2}}$.

\smallskip

$\mathcal{RP}(\mathbb{F}_{3})$ is a cyclic $R_{\mathbb{F}_{3}}$-module generated by $[-1]$ subject to the relation $[1]=0$ and

\begin{center}
$0=2\psi_{1}(-1)=2([-1]+\langle-1\rangle[-1])$.
\end{center}

$\mathcal{RB}(\mathbb{F}_{3})=\mathcal{RP}_{1}(\mathbb{F}_{3})$ is the submodule generated by $\psi_{1}(-1)=[-1]+\langle-1\rangle[-1]$.

\smallskip

For any field $F$, there is a natural surjective homomorphism of $R_{F}$-modules $H_{3}(\mathrm{SL}_{2}(F),\mathbb{Z})\rightarrow K_{3}^{\mathrm{ind}}(F)$, it can be shown this is an homomorphism of $R_{F}$-modules where $F^{\times}/(F^{\times})^{2}$ acts trivially on $K_{3}^{\mathrm{ind}}(F)$. In \cite{mirzaii:third}, we get this map induces an isomorphism

\begin{center}
$H_{0}\left(F^{\times}/(F^{\times})^{2},H_{3}\left(\mathrm{SL}_{2}(F),\mathbb{Z}\left[\frac{1}{2}\right]\right)\right)\cong K_{3}^{\mathrm{ind}}(F)\left[\frac{1}{2}\right]$
\end{center}

Now let $H_{3}(\mathrm{SL}_{2}(F),\mathbb{Z})_{0}$ denote the kernel of the surjective homomorphism $H_{3}(\mathrm{SL}_{2}(F),\mathbb{Z})\rightarrow K^{\mathrm{ind}}_{3}(F)$. Note the isomorphism above implies

\begin{center}
$H_{3}\left(\mathrm{SL}_{2}(F),\mathbb{Z}\left[\frac{1}{2}\right]\right)_{0}=\mathcal{I}_{F}H_{3}\left(\mathrm{SL}_{2}(F),\mathbb{Z}\left[\frac{1}{2}\right]\right)$
\end{center}

\begin{Rec}\label{Recall_Number_fields_K3_ind}
Let us recall that when $F$ is a number field the surjective homomorphism $H_{3}(\mathrm{SL}_{2}(F),\mathbb{Z})\rightarrow K^{\mathrm{ind}}_{3}(F)$ is split as a map of abelian groups. In fact $K_{3}^{ind}(F)$ is finitely generated abelian group an it is enough that there is a torsion subgroup of $H_{3}(\mathrm{SL}_{2}(F),\mathbb{Z})$ mapping isomorphically to the cyclic torsion subgroup of $K_{3}^{ind}(F)$. This latter statement follows from the explicit calculations of Christian Zickert in \cite[Section 8]{CZickert} . It follows as an abelian group

\begin{center}
$H_{3}(\mathrm{SL}_{2}(F),\mathbb{Z})\cong K_{3}^{ind}(F)\oplus H_{3}(\mathrm{SL}_{2}(F),\mathbb{Z})_{0}$
\end{center}

for any number field $F$.
\end{Rec}

Furthermore, we have the following proposition:

\begin{prop}\cite[Corollary 2.8, Corollary 4.4]{Articlediscretevaluation}\label{PropsitionRP+(F)andH3(SL2)}
For any field $F$, the map $H_{3}(\mathrm{SL}_{2}(F),\mathbb{Z})\rightarrow \mathcal{RP}(F)$ induces an isomorphism of $R_{F}$-modules

\begin{center}
$H_{3}\left(\mathrm{SL}_{2}(F),\mathbb{Z}\left[\frac{1}{2}\right]\right)_{0}=\mathcal{I}_{F}H_{3}\left(\mathrm{SL}_{2}(F),\mathbb{Z}\left[\frac{1}{2}\right]\right)\cong\mathcal{I}_{F}\mathcal{RP}_{1}(F)\left[\frac{1}{2}\right]$
\end{center}

and furthermore

\begin{center}
$\mathcal{RP}_{1}(F)\left[\frac{1}{2}\right]=\widetilde{\mathcal{RP}}_{1}(F)\left[\frac{1}{2}\right]=e_{+}^{-1}\widetilde{\mathcal{RP}}(F)\left[\frac{1}{2}\right]$
\end{center}

where $e_{+}^{-1}:=\frac{p_{+}^{-1}}{2}=\frac{1+\langle-1\rangle}{2}\in R_{F}\left[\frac{1}{2}\right]$
\end{prop}

\begin{flushright}
    $\Box$
\end{flushright}

Note that from $\mathcal{RB}(F)\subseteq \mathcal{RP}_{1}(F)$ and Theorem \ref{TheoremBlochWignerComplex}, it follows that the square class $\langle-1\rangle$ acts trivially on $H_{3}\left(\mathrm{SL}_{2}(F),\mathbb{Z}\left[\frac{1}{2}\right]\right)$ (see \cite[Corollary 4.6]{Articlediscretevaluation} for more details). Now we define $\mathcal{RP}_{+}(F)$ to be $\widetilde{\mathcal{RP}}(F)$ modulo the submodule generated by the elements $(1-\langle-1\rangle)[x]$, $x\in F^{\times}$. Thus $\mathcal{RP}_{+}(F)$ is the $R_{F}$-module generated by the elements $[x]$, $x\in F^{\times}$ subject to the relations

\begin{itemize}
\item[1)] $[1]=0$.
\item[2)] $S_{x,y}=0$ for $x,y\not=1$.
\item[3)] $\langle-1\rangle[x]=[x]$ for all $x\in F^{\times}$.
\item[4)] $[x]=-[x^{-1}]$ for all $x\in F^{\times}$.
\end{itemize}

The proposition above implies that the map $H_{3}(\mathrm{SL}_{2}(F),\mathbb{Z})\rightarrow \mathcal{RP}(F)$ induces an isomorphism

\begin{center}
$H_{3}\left(\mathrm{SL}_{2}(F),\mathbb{Z}\left[\frac{1}{2}\right]\right)_{0}\cong\mathcal{I}_{F}\mathcal{RP}_{+}(F)\left[\frac{1}{2}\right]$
\end{center}

\begin{Rem}
It will be convinient below to introduce the following notation in $\mathcal{RP}_{+}(F)$:

\begin{center}
$[1]:=0$,\;\;\;\;$[0]:=C_{F}$ \;\; and\;\; $[\infty]:=-C_{F}$.
\end{center}

Thus the symbol $[x]\in \mathcal{RP}_{+}(F)$ is defined for all $x\in \mathbb{P}^{1}(F)$. Furthermore $[x]+[1-x]=C_{F}$ and $[x]=-[x^{-1}]$ in $\mathcal{RP}_{+}(F)$ for all $x\in \mathbb{P}^{1}(F)$.
\end{Rem}

\section{Characters}\label{Section_Characters}

In this section we will review the character-theoretic method introduced by Hutchinson (\cite{Articlediscretevaluation}) for calculating with
modules over a group ring with an elementary abelian $2$-group. Furthermore,
we prove some new related results about subrings in number field which will be useful in Section \ref{Section_proof_of_the_main_theorem}.

\bigskip

Let $G$ be an abelian group satisfying $g^{2}=1$ for all $g\in G$. Let
$\mathcal{R}$ denote the group ring $\mathbb{Z}[G]$. For a character $\chi\in
\widehat{G}:=\mathrm{Hom}(G,\mu_{2})$, let $\mathcal{R}^{\chi}$ be the ideal
of $\mathcal{R}$ generated by the elements $\{g-\chi(g)\;|\; g\in G\}$; i.e.
$\mathcal{R}^{\chi}$ is the kernel of the ring homomorphism
$\rho(\chi):\mathcal{R}\rightarrow \mathbb{Z}$ sending $g$ to $\chi(g)$ for
any $g\in G$. We let $\mathcal{R}_{\chi}$ denote the associated
$\mathcal{R}$-algebra structure on $\mathbb{Z}$. In other words
$\mathcal{R}_{\chi}:=\mathcal{R}/\mathcal{R}^{\chi}$.

\bigskip

If $M$ is a $\mathcal{R}$-module, we let $M^{\chi}=\mathcal{R}^{\chi}M$ and
we let

\begin{center}
$M_{\chi}:=M/M^{\chi}=(\mathcal{R}/\mathcal{R}^{\chi})\otimes_{\mathcal{R}}M=\mathcal{R}_{\chi}\otimes_{\mathcal{R}}M$.
\end{center}

Thus $M_{\chi}$ is the largest quotient module of $M$ with the property that
$g\cdot m=\chi(g)\cdot m$ for all $g\in G$. In particular, if $\chi=\chi_{0}$
is the trivial character, then $\mathcal{R}^{\chi_{0}}$ is the augmentation
ideal $\mathcal{I}_{G}$, $M^{\chi_{0}}=\mathcal{I}_{G}M$ and
$M_{\chi_{0}}=M_{G}$.

\bigskip

Given $m\in M$, $\chi\in \widehat{G}$, we denote the image of $m$ in
$M_{\chi}$ by $m_{\chi}$. For example for any character $\chi\in
\widehat{F^{\times}/(F^{\times})^{2}}$, we can give a presentation of the
$R_{F}$-module $\mathcal{RP}_{+}(F)_{\chi}$ which is our main object of
study. $\mathcal{RP}_{+}(F)_{\chi}$ is the $R_{F}$-module with generators
$[x]_{\chi}$, $x\in F^{\times}$, subject to the relations

\begin{itemize}
\item[1)] $\langle a\rangle\cdot[x]_{\chi}:=\chi(a)\cdot[x]_{\chi}$ for
    all $a,x\in F^{\chi}$.
\item[2)] $[1]_{\chi}=0$
\item[3)] The five term relation

\begin{center}
$0=[x]_{\chi}-[y]_{\chi}+\chi(x)\left[\frac{y}{x}\right]_{\chi}-\chi(x^{-1}-1)\left[\frac{1-x^{-1}}{1-y^{-1}}\right]_{\chi}+\chi(1-x)\left[\frac{1-x}{1-y}\right]_{\chi}$
\end{center}

for all $x,y\not=1$

\item[4)] $\chi(-1)\cdot[x]_{\chi}=[x]_{\chi}$ for all $x\in F^{\times}$.
\item[5)] $[x]_{\chi}=-[x^{-1}]_{\chi}$ for all $x\in F^{\times}$.
\end{itemize}

\begin{Rec}
Let us recall for $M$ a $\mathcal{R}$-module, $M_{\chi}=\mathcal{R}_{\chi}\otimes_{\mathcal{R}}M$. Hence, for any
$\mathcal{R}$-homomorphism $f:M\rightarrow N$, it follows that we have the
following induced map

\begin{center}
$f_{\chi}=Id_{\mathcal{R}_{\chi}}\otimes
f:M_{\chi}=\mathcal{R}_{\chi}\otimes_{\mathcal{R}}M\rightarrow
\mathcal{R}_{\chi}\otimes_{\mathcal{R}}N=N_{\chi}$
\end{center}

\end{Rec}

The following the character-theoretic local-global principle is central to
calculations in the following chapter below

\begin{prop}\cite[Section 3]{Articlediscretevaluation}\label{Proposition_Local_global_principle_character}
\begin{itemize}
\item[(1)] For any $\chi\in\widehat{G}$, $M\rightarrow M_{\chi}$ is an
    exact functor on the category of
    $\mathcal{R}\left[\frac{1}{2}\right]$-modules.
\item[(2)] Let $f:M\rightarrow N$ be an
    $\mathcal{R}\left[\frac{1}{2}\right]$-module homomorphism. For any
    $\chi\in\widehat{G}$, let $f_{\chi}:M_{\chi}\rightarrow N_{\chi}$ be
    the induced map. Then $f$ is bijective (resp. injective, surjective)
    if and only if $f_{\chi}$ is bijective (resp. injective, surjective)
    for all $\chi\in\widehat{G}$.
\end{itemize}
\end{prop}

\begin{flushright}
    $\Box$
\end{flushright}


\begin{lema}\cite[Lemma 5.3]{ArticlehomologyQprint}
Let $F$ be a field. Let $\chi\in \widehat{F^{\times}/(F^{\times})^{2}}$.
Suppose that $a\in F^{\times}$ satisfy $\chi(1-a)=-1$ and $\chi(a)=1$ then
$[a]_{\chi}=0$ in $\mathcal{RP}_{+}(F)\left[\frac{1}{2}\right]_{\chi}$
\end{lema}

\begin{flushright}
    $\Box$
\end{flushright}




\begin{cor}\cite[Lemma 5.4, Corollary 5.5]{ArticlehomologyQprint}\label{ArticlehomologyQprint(1-l)a}
Let $F$ be a field. Let $\chi\in \widehat{F^{\times}/(F^{\times})^{2}}$ with
$\chi(-1)=1$. Suppose that $\ell\in F^{\times}$ satisfy $\chi(\ell)=-1$ and
$\chi(1-\ell)=1$. Then $[a]_{\chi}=[(1-\ell)^{m}a]_{\chi}$ in
$\mathcal{RP}_{+}(F)\left[\frac{1}{2}\right]_{\chi}$ for all $a\in
\mathbb{P}^{1}(F)$ and all $m\in\mathbb{Z}$.
\end{cor}

\begin{flushright}
    $\Box$
\end{flushright}


\begin{Rem}\label{Remark[a]=[la]}
For a field $F$, let $\chi\in \widehat{F^{\times}/(F^{\times})^{2}}$ such
that $\chi(-1)=1$. Suppose that $\ell\in F^{\times}$ satisfy $\chi(\ell)=-1$
and $\chi(1-\ell)=-1$ then $\chi\left(\frac{1}{\ell}\right)=-1$,
$\chi\left(1-\frac{1}{\ell}\right)=1$. Therefore by Corollary
\ref{ArticlehomologyQprint(1-l)a}, we can deduce the following for all $a\in
\mathbb{P}^{1}(F)$:

\begin{itemize}
\item[a)] If $x(\ell)=1$ and $\chi(1-\ell)=-1$ then $[a]_{\chi}=[\ell
    a]_{\chi}$.
\item[b)] If $\chi(\ell)=-1$ and $\chi(1-\ell)=-1$ then
    $[a]_{\chi}=\left[\left(1-\frac{1}{\ell}\right)a\right]_{\chi}$.
\end{itemize}
\end{Rem}

From Corollary \ref{ArticlehomologyQprint(1-l)a}, we can deduce the following
result:

\begin{cor}\cite[Lemma 5.6]{ArticlehomologyQprint}\label{CorollarySumQprint}
Let $F$ be a field with at least $4$ elements. Let $\chi\in \widehat{F^{\times}/(F^{\times})^{2}}$ with
$\chi(-1)=1$. Suppose that $\ell\in F^{\times}$ satisfy $\chi(\ell)=-1$ and
$\chi(1-\ell)=1$. Then $[a]_{\chi}=[a+t\ell]_{\chi}$ in
$\mathcal{RP}_{+}(F)\left[\frac{1}{2}\right]_{\chi}$ for all $a\in
F^{\times}$ and all $t\in\mathbb{Z}$.
\end{cor}

\begin{flushright}
    $\Box$
\end{flushright}

Now let $R$ be a commutative ring and $x\in R$. Let us recall that $\mathbb{Z}[x]$ is
the subring

\begin{center}
$\{a_{0}+a_{1}x+\cdots+a_{n}x^{n}\;|\; n\geq0\;\; a_{i}\in\mathbb{Z}\}\subset
R$.
\end{center}

Now $x\mathbb{Z}[x]=\{xp(x)\;|\; p(x)\in \mathbb{Z}[x] \}$; i.e.
$x\mathbb{Z}[x]$ is the ideal in $\mathbb{Z}[x]\subset R$ generated by $x\in
R$.

\bigskip

The following proposition is a generalization of the Corollary
\ref{CorollarySumQprint} for any field $F$:

\begin{prop}\label{Propositionelement_lZ[l]}
Let $F$ be a field with at least $4$ elements. Let $\chi\in \widehat{F^{\times}/(F^{\times})^{2}}$ with
$\chi(-1)=1$. If $\chi(\ell)=-1$, $\chi(1-\ell)=1$ for $\ell\not=1\in F$. Then $[a]_{\chi}=[a+t]_{\chi}$ in $\mathcal{RP}_{+}(F)\left[\frac{1}{2}\right]_{\chi}$ for all $a\in F^{\times}$ and all $t\in \ell\mathbb{Z}[\ell]$ 
\end{prop}

$\textbf{\emph{Proof}}.$ Let $t\in \ell\mathbb{Z}[\ell]$. So there exist
$m\geq1$, $b_{1},\ldots, b_{m}\in\mathbb{Z}$ with
$t=b_{1}\ell+\cdots+b_{m}\ell^{m}$.  We will prove it by induction on $m$. If
$m=1$, the result follows from the Corollary \ref{CorollarySumQprint}. Now we
suppose it is valid for any number $m\leq k$. We will prove it for the case
$k+1$. First, note that $\chi\left(\frac{1}{1-\ell}\right)=1$ and
$\chi\left(1-\frac{1}{1-\ell}\right)=\chi(-1)\chi(\ell)\chi\left(\frac{1}{1-\ell}\right)=-1$.
Thus

\begin{eqnarray}
\nonumber [a]_{\chi} &= &\left[\frac{a}{1-\ell}\right]_{\chi}  \mbox{(From the Remark \ref{Remark[a]=[la]})}\\
\nonumber &=&\left[\frac{a}{1-\ell}+f_{1}\ell+\cdots +f_{k}\ell^{k}\right]_{\chi} \mbox{ (For any $f_{1},\ldots f_{k}\in\mathbb{Z}$ by the induction hypothesis)}\\
\nonumber &=&\left[a+(1-\ell)(f_{1}\ell+\cdots +f_{k}\ell^{k})\right]_{\chi}\mbox{(From the Corollary \ref{ArticlehomologyQprint(1-l)a})}\\
\nonumber &=&\left[a+(1-\ell)(f_{1}\ell+\cdots +f_{k}\ell^{k})+   c_{1}\ell+\cdots +c_{k}\ell^{k} \right]_{\chi}\mbox{ (For any $c_{1},\ldots c_{k}\in\mathbb{Z}$ by hypothesis)}\\
\nonumber &=&\left[a+ (f_{1}+c_{1})\ell+ (f_{2}-f_{1}+c_{2})\ell^{2}+\cdots +(f_{k}-f_{k-1}+c_{k})\ell^{k-1}+f_{k}\ell^{k+1} \right]_{\chi}.
\end{eqnarray}

Since $f_{1},\ldots f_{k}, c_{1},\ldots, c_{k}$ are chosen arbitrarily
integers. It follows that for any $b_{1},\ldots, b_{k+1}\in\mathbb{Z}$, we can solve the following system of equations:

\begin{eqnarray}
\nonumber f_{1}+c_{1} &=&b_{1}\\
\nonumber f_{2}-f_{1}+c_{2}&=&b_{2}\\
\nonumber &\vdots&\\
\nonumber f_{k}-f_{k-1}+c_{k}&=&b_{k}\\
\nonumber f_{k}&=&b_{k+1}
\end{eqnarray}

Therefore $[a]_{\chi}=\left[a+b_{1}\ell+\cdots
+b_{k+1}\ell^{k+1}\right]_{\chi}$ for any $b_{1},\ldots,
b_{k+1}\in\mathbb{Z}$.

\begin{flushright}
    $\Box$
\end{flushright}

\begin{lema}\label{Corollary_lZ[l]_NumberField}
Let $F$ be a number field such that $[F:\mathbb{Q}]=d$ and $\mathcal{O}_{F}$
its ring of algebraic integers. Let $\ell\in \mathcal{O}_{F}$ of degree $d$
and $N=N_{F/\mathbb{Q}}(\ell)\in \mathbb{Z}$. Then

\begin{center}
$\ell\mathbb{Z}[\ell]=\mathbb{Z}N+\mathbb{Z}\ell+\cdots\mathbb{Z}\ell^{d-1}\subseteq
\mathcal{O}_{F}$.
\end{center}

\end{lema}

$\textbf{\emph{Proof}}.$  Let $p(t)=N+b_{1}t+\cdots + b_{d-1}t^{d-1}+t^{d}$
the minimal polynomial of $\ell$. Then

\begin{center}
$0=p(\ell)=N+b_{1}\ell+\cdots +b_{d-1}\ell^{d-1}+\ell^{d}$. 
\end{center}

Since $b_{1}\ell+\cdots b_{d-1}\ell^{d-1}+\ell^{d}\in \ell\mathbb{Z}[\ell]$,
it follows that $N\in \ell\mathbb{Z}[\ell]$. Therefore

\begin{center}
$\mathbb{Z}N+\mathbb{Z}\ell+\cdots\mathbb{Z}\ell^{d-1}\subseteq
\ell\mathbb{Z}[\ell]$.
\end{center}

Now note that $\ell^{k}\in
\mathbb{Z}N+\mathbb{Z}\ell+\cdots\mathbb{Z}\ell^{d-1}$ with
$k\in\{1,2,\ldots,d-1\}$ and from the equation $0=p(\ell)$ we get that
$\ell^{d}\in \mathbb{Z}N+\mathbb{Z}\ell+\cdots\mathbb{Z}\ell^{d-1}$.

\bigskip

Now, multiplying the equation $0=p(\ell)$ by $\ell$ we obtain that

\begin{center}
$\ell^{d+1}=-N-b_{1}\ell-\cdots -b_{d-1}\ell^{d}$.
\end{center}

Since $\ell^{d}\in \mathbb{Z}N+\mathbb{Z}\ell+\cdots\mathbb{Z}\ell^{d-1}$ ,
it follows that
$\ell^{d+1}\in\mathbb{Z}N+\mathbb{Z}\ell+\cdots\mathbb{Z}\ell^{d-1}\subseteq
\ell\mathbb{Z}[\ell]$. Hence following a similar recursive way, we get that
$\ell^{m}\in\mathbb{Z}N+\mathbb{Z}\ell+\cdots\mathbb{Z}\ell^{d-1}\subseteq
\ell\mathbb{Z}[\ell]$ with $m\geq d+1$. Therefore it follows that
$\ell\mathbb{Z}[\ell]
\subseteq\mathbb{Z}N+\mathbb{Z}\ell+\cdots\mathbb{Z}\ell^{d-1}$.

\begin{flushright}
    $\Box$
\end{flushright}

\begin{lema}\label{corollary_1_isintheset}
Let $F$ be a number field such that $[F:\mathbb{Q}]=d$ and $\mathcal{O}_{F}$
its ring of algebraic integers. Let $\ell\in \mathcal{O}_{F}$ of degree $d$
and $N=N_{F/\mathbb{Q}}(\ell)\in \mathbb{Z}$. Then $\mathbb{Z}[\ell]
\subset\frac{1}{\ell}\mathbb{Z}\left[\frac{1}{\ell}\right]$.
\end{lema}

$\textbf{\emph{Proof}}.$  Let $p(t)=N+b_{1}t+\cdots b_{d-1}t^{d-1}+t^{d}$ the
minimal polynomial of $\ell$. Then

\begin{center}
$0=p(\ell)=N+b_{1}\ell+\cdots b_{d-1}\ell^{d-1}+\ell^{d}$.
\end{center}

Multiplying the last equation by $\ell^{-d}$, we get

\begin{center}
$0=p(\ell)=\frac{N}{\ell^{d}}+\frac{b_{1}}{\ell^{d-1}}+\cdots
\frac{b_{d-1}}{\ell}+1$.
\end{center}

Since $\frac{N}{\ell^{d}}+\frac{b_{1}}{\ell^{d-1}}+\cdots
\frac{b_{d-1}}{\ell}\in \frac{1}{\ell}\mathbb{Z}\left[\frac{1}{\ell}\right]$.
It follows that $1\in \frac{1}{\ell}\mathbb{Z}\left[\frac{1}{\ell}\right]$
whence we obtain that $\mathbb{Z}\subset
\frac{1}{\ell}\mathbb{Z}\left[\frac{1}{\ell}\right]$.

\bigskip

 Now we can finally prove $\mathbb{Z}[\ell] \subset\frac{1}{\ell}\mathbb{Z}\left[\frac{1}{\ell}\right]$. For that, it is enough to prove that $\ell^{m}\in \frac{1}{\ell}\mathbb{Z}\left[\frac{1}{\ell}\right]$ where $m\geq1$. We will prove it by induction on $m$. If $m=1$, note that since $1\in \frac{1}{\ell}\mathbb{Z}\left[\frac{1}{\ell}\right]$, we get that $1=\frac{a_{1}}{\ell}+\frac{a_{2}}{\ell^{2}}+\cdots +\frac{a_{d-1}}{\ell^{d-1}}$ with $a_{1},a_{2},\ldots,a_{d-1}\in\mathbb{Z}$, then multiplying the last equation by $\ell$ we get the following equation:

\begin{center}
$\ell=a_{1}+\frac{a_{2}}{\ell}+\cdots +
\frac{a_{d-2}}{\ell^{d-3}}+\frac{a_{d-1}}{\ell^{d-2}}$.
\end{center}

Since $a_{1}\in \frac{1}{\ell}\mathbb{Z}\left[\frac{1}{\ell}\right]$, it
follows that $\ell\in \frac{1}{\ell}\mathbb{Z}\left[\frac{1}{\ell}\right]$.
Now we suppose that $\ell^{m}\in \mathbb{Z}\left[\frac{1}{\ell}\right]$ for all $m\leq k$. We will prove the case
$m=k+1$. Multiplying the last equation by $\ell^{k}$, we get the equation

\begin{center}
$\ell^{k+1}=a_{1}\ell^{k}+a_{2}\ell^{k-1}\cdots + a_{d-1}\ell^{k-d+2}$.
\end{center}

By the induction hypothesis, it follows that $\ell^{k+1}\in
\frac{1}{\ell}\mathbb{Z}\left[\frac{1}{\ell}\right]$. Therefore as required  $\mathbb{Z}[\ell]
\subset\frac{1}{\ell}\mathbb{Z}\left[\frac{1}{\ell}\right]$.

\begin{flushright}
    $\Box$
\end{flushright}

\begin{Rem}\label{Remark_for_corollary_l_isintheset}
From the proof of Lemma \ref{corollary_1_isintheset}, we have that
$\frac{1}{\ell}\mathbb{Z}\left[\frac{1}{\ell}\right]$ is closed under
multiplication by $\ell$.
\end{Rem}

\section{Fields with a valuation}\label{Section_Fields_with_a_valuation}

Given a field $F$ and a surjective valuation $\upsilon:F^{\times}\rightarrow \Gamma$, where $\Gamma$ is a totally ordered additive abelian group, we let $\mathcal{O}_{\upsilon}:=\{x\in F^{\times}\;|\;\upsilon(x)\geq0\}\cup\{0\}$ be the associated valuation ring, with maximal ideal $\mathcal{M}_{\upsilon}=\{x\in \mathcal{O}_{\upsilon}\;|\;\upsilon(x)\not=0\}$, group of units $U_{\upsilon}=U:=\mathcal{O}_{\upsilon}\backslash\mathcal{M}_{\upsilon}$ and residue field $k=k(\upsilon):=\mathcal{O}_{\upsilon}/\mathcal{M}_{\upsilon}$.

\smallskip

Since $\Gamma$ is a torsion-free group, we have a short exact sequence of $\mathbb{F}_{2}$-vector spaces

\begin{center}
$\xymatrix{1\ar[r]&\frac{U}{U^{2}}\ar[r]&\frac{F^{\times}}{(F^{\times})^{2}}\ar[r]&\frac{\Gamma}{2}\ar[r]&1}$.
\end{center}

We have a homomorphism of commutative rings

\begin{center}
$\xymatrix{\mathbb{Z}\left[\frac{U}{U^{2}}\right]\ar@{^{(}->}[r] \ar@{>>}[d]&R_{F}\\
R_{k}&}$
\end{center}

\smallskip

Therefore given an $R_{k}$-module $M$, we denote by $\mathrm{Ind}_{k}^{F}M$ the $R_{F}$-module $R_{F}\otimes_{\mathbb{Z}[U_{\upsilon}]}M$. 
We have the following result \cite[Section 5]{Articlediscretevaluation}:

\begin{lema}
There is a natural homomorphism of $R_{F}$-modules $S_{\upsilon}:\widetilde{\mathcal{RP}}(F)\rightarrow \mathrm{Ind}_{k}^{F}\widetilde{\mathcal{RP}}(k)$ given by

\begin{center}
$S_{\upsilon}([x])=\left\{
                     \begin{array}{ll}
                       1\otimes[\overline{x}], & \hbox{$\upsilon(x)=0$;} \\
                       1\otimes C_{k}, & \hbox{$\upsilon(x)>0$;} \\
                       -(1\otimes C_{k}), & \hbox{$\upsilon(x)<0$.}
                     \end{array}
                   \right.
$
\end{center}

\end{lema}

\begin{flushright}
    $\Box$
\end{flushright}












\subsection{Discrete valuations}

Suppose that $\upsilon:F^{\times}\rightarrow \mathbb{Z}$ is a discrete valuation on the field $F$ with residue field $k=k(\upsilon)$. Let $\chi_{\upsilon}:F^{\times}/(F^{\times})^{2}\rightarrow \mu_{2}$ denote the associated character defined by $\chi_{\upsilon}(a)=(-1)^{\upsilon(a)}$. For an abelian group $M$, we let $M\{\upsilon\}$ denote the $R_{F}$-module $R_{\chi_{\upsilon}}\otimes_{\mathbb{Z}}M$. Equivalently we equip $M$ with the $R_{F}$-module structure $\langle a\rangle m:=(-1)^{\upsilon(a)}m$ for all $a\in F^{\times}$ and $m\in M$. We have the following result:

\begin{prop}\cite[Theorem 3.7]{ArticlehomologyQprint}\label{NaturalIsomorphismRP+(F)}
Let $F$ be a field with discrete valuation $\upsilon:F^{\times}\rightarrow \mathbb{Z}$ and residue field $k$. Then we have a natural isomorphism

\begin{center}
$\xymatrix{\mathcal{RP}_{+}(F)\left[\frac{1}{2}\right]_{\chi_{\upsilon}}\ar[r]^{S_{\upsilon}}&\left(\mathrm{Ind}^{F}_{k}\mathcal{RP}_{+}(k)\left[\frac{1}{2}\right]\right)_{\chi_{\upsilon}}\ar[r]^{\cong}& \mathcal{P}(k)\left[\frac{1}{2}\right]\{\upsilon\} }$
\end{center}

\end{prop}

\begin{flushright}
    $\Box$
\end{flushright}

For any field $F$ with discrete valuation $\upsilon$, we let $\overline{S}_{\upsilon}$ denote the composite $R_{F}$-module homomorphism

\begin{center}
$H_{3}\left(\mathrm{SL}_{2}(F), \mathbb{Z}\left[\frac{1}{2}\right]\right)\rightarrow \left(\mathcal{RP}_{+}(F)\left[\frac{1}{2}\right]\right)_{\chi_{\upsilon}}\cong\mathcal{P}(k(\upsilon))\left[\frac{1}{2}\right]\{\upsilon\}$
\end{center}

By abuse of notation, we will use the same symbol to denote $\overline{S}_{\upsilon}$ restricted to $H_{3}\left(\mathrm{SL}_{2}(F), \mathbb{Z}\left[\frac{1}{2}\right]\right)_{0}$.

\begin{Rem}\label{squareclassuniformizer}
Let us recall that if $\pi$ is a uniformizer for the valuation $\upsilon$,
then $\upsilon(\pi)=1$. Therefore the square class of $\pi$ acts as $-1$ on
the factor $\mathcal{P}(k(\upsilon))\{\upsilon\}$ on the right.
\end{Rem}

\begin{Rec}
Give a family $\mathcal{V}$ of discrete valuations of $F$, we obtain a map

\begin{center}
$\xymatrix{H_{3}(\mathrm{SL}_{2}(F),\mathbb{Z} )\ar[r]& \mathcal{RP}_{+}(F)\ar[r]& \displaystyle\prod_{\upsilon\in \mathcal{V}}\mathcal{P}(k(\upsilon))\{\upsilon\}}$.
\end{center}

Note that when we restrict to $H_{3}(\mathrm{SL}_{2}(F),\mathbb{Z})_{0}$ and tensor with $\mathbb{Z}\left[\frac{1}{2}\right]$, the image lies in the direct sum instead (For more details see \cite[Theorem 5.1]{ArticleBloch} ).

\end{Rec}

Let $\mathcal{V}$ be a set of discrete valuations, we get a homomorphism $\overline{S}=\{\overline{S}_{\upsilon}\}_{\upsilon\in\mathcal{V}}$  of $R_{F}$-modules

\begin{center}
$\overline{S}:H_{3}\left(\mathrm{SL}_{2}(F),\mathbb{Z}\left[\frac{1}{2}\right]\right)_{0}\rightarrow \displaystyle\bigoplus_{v\in\mathcal{V}}\mathcal{P}(k(\upsilon))\{\upsilon\}$.
\end{center}

If $\upsilon(x)$ is odd, then $\langle x\rangle$ acts as $-1$ on $M\{\upsilon\}$.


\section{Quadratic fields.}\label{Section_Quadratic_fields}

Let us recall that if $a,b\in R$ then $a$ is associated to $b$ ($a\sim b$), if there exist $u\in R^{\times}$ such
that $a=ub$.



For the field $\mathbb{Q}(\sqrt{m})$, we have the ring of algebraic
integers
$\mathcal{O}_{\mathbb{Q}(\sqrt{m})}=\mathbb{Z}[\omega_{m}]=\{a+\omega_{m}b\;|\;
a,b\in\mathbb{Z}\}$ where

\begin{center}
$\omega_{m}=\left\{
              \begin{array}{ll}
                \sqrt{m}, & \hbox{$m\;\equiv\;2,3\;(\mathrm{mod}\;4)$;} \\
                \frac{1+\sqrt{m}}{2}, & \hbox{$m\;\equiv\;1\;(\mathrm{mod}\;4)$.}
              \end{array}
            \right.
$
\end{center}



For $\alpha\in \mathbb{Q}(\sqrt{m})$, let $\overline{\alpha}$ denote its conjugate. Then if $\alpha=a+\omega_{m}b\in \mathcal{O}_{\mathbb{Q}(\sqrt{m})}$ , we get that:

\begin{center}
$\overline{\alpha}=\left\{
                                                          \begin{array}{ll}
                                                            a-\omega_{m}b, & \hbox{$m\;\equiv\;2,3\;(\mathrm{mod}\;4)$;} \\
                                                            (a+b)-\omega_{m}b, & \hbox{$m\;\equiv\;1\;(\mathrm{mod}\;4)$.}
                                                          \end{array}
                                                        \right.$
\end{center}

Therefore in any quadratic field $\mathbb{Q}(\sqrt{m})$, the norm of
$\alpha\in\mathcal{O}_{\mathbb{Q}(\sqrt{m})}$ is given by

\begin{center}
$N_{\mathbb{Q}(\sqrt{m})/\mathbb{Q}}(\alpha)=N(\alpha):=\alpha\overline{\alpha}=\left\{
                                                          \begin{array}{ll}
                                                            a^{2}-mb^{2}, & \hbox{$m\;\equiv\;2,3\;(\mathrm{mod}\;4)$;} \\
                                                            a^{2}+ab+\left(\frac{1-m}{4}\right)b^{2}, & \hbox{$m\;\equiv\;1\;(\mathrm{mod}\;4)$.}
                                                          \end{array}
                                                        \right.
$.
\end{center}

Let us recall that the \emph{discriminant} $\delta$ of any quadratic field $\mathbb{Q}(\sqrt{m})$
is:

\begin{center}
$\delta=\left\{
          \begin{array}{ll}
            m, & \hbox{$m\;\equiv\;1\;(\mathrm{mod}\;4)$;} \\
            4m, & \hbox{$otherwise$.}
          \end{array}
        \right.
$
\end{center}

Now, the prime elements in $\mathcal{O}_{\mathbb{Q}(\sqrt{m})}$ are found by
implementing the three following primality criteria 

\begin{itemize}
\item If $N(\alpha)$ is a rational prime, then $\alpha$ is a prime.
\item An odd rational prime $p$ is a prime in
    $\mathcal{O}_{\mathbb{Q}(\sqrt{m})}$ if and only if $p$ does not
    divide $\delta$ and $\left(\frac{\delta}{p}\right)=-1$.
\item $2$ is prime in $\mathcal{O}_{\mathbb{Q}(\sqrt{m})}$ if and only if
    $m\;\equiv\;5\;(\mathrm{mod}\;8)$.
\end{itemize}

We introduce the following notation. For $\alpha=a+b\omega_{m}\in
\mathbb{Z}[\omega_{m}]$, we let $a:=R(\alpha)$ and $b:=I(\alpha)$. Then
with this notation we have:

\begin{itemize}
\item[i)] If $m\;\equiv\;2,3\;(\mathrm{mod}\;4)$ then
    \begin{center}
    $R(\overline{\alpha})=R(\alpha)$\;\;\; and\;\;\;
    $I(\overline{\alpha})=-I(\alpha)$.
    \end{center}
\item[ii)]  If $m\;\equiv\;1\;(\mathrm{mod}\;4)$ then

    \begin{center}
    $R(\overline{\alpha})=R(\alpha)+I(\alpha)$\;\;\; and\;\;\;
    $I(\overline{\alpha})=-I(\alpha)$.
    \end{center}

\end{itemize}

Thus

\begin{center}
$\overline{\alpha}=\left\{
                                                          \begin{array}{ll}
                                                            R(\alpha)-\omega_{m}I(\alpha), & \hbox{$m\;\equiv\;2,3\;(\mathrm{mod}\;4)$;} \\
                                                            (R(\alpha)+I(\alpha))-\omega_{m}I(\alpha), & \hbox{$m\;\equiv\;1\;(\mathrm{mod}\;4)$.}
                                                          \end{array}
                                                        \right.$
\end{center}

We will require the following observation below.





\begin{lema}\label{R(alpha)_and__N(alpha)_are_not_coprimes}
Let $\mathbb{Q}(\sqrt{-m})$ with $m>0$ such that $m\;\equiv\;1,2\;(\mathrm{mod}\;4)$ and let
$\alpha=R(\alpha)+\omega_{m}I(\alpha)\in \mathcal{O}_{\mathbb{Q}(\sqrt{-m})}$ a prime
element such that $N(\alpha)=p$ is a rational prime. If $p$ divides
$R(\alpha)$ or $m$ in $\mathbb{Z}$ then $m=p$. In particular $\alpha\sim\sqrt{-p}$.
\end{lema}

$\textbf{\emph{Proof}}.$ If $p|R(\alpha)$, we have $R(\alpha)=pn$ for some
$n\in\mathbb{Z}$. Thus

\begin{center}
$p=N(\alpha)=p^{2}n^{2}+m(I(\alpha))^{2}$.
\end{center}

From this the last equation, we get $m=p$ and $R(\alpha)=0$,
$I(\alpha)=\pm1$ as otherwise $N(\alpha)=(R(\alpha))^{2}+m(I(\alpha))^{2}>p$.

\smallskip

If $p|m$ then $m=p\ell$  for some $\ell\in\mathbb{Z}$. Thus $p=N(\alpha)=(R(\alpha))^{2}+p\ell(I(\alpha))^{2}$. It follows that $p|R(\alpha)$ and from above, we get our result.

\begin{flushright}
    $\Box$
\end{flushright}

\section{The main theorem}\label{Section_proof_of_the_main_theorem}

Let $F$ be a number field. We consider $\Pi$, the set of discrete valuations on $F$. As we seen, the maps $S_{p}$ induce



\begin{center}
$\overline{S}:H_{3}\left(\mathrm{SL}_{2}(F),\mathbb{Z}\left[\frac{1}{2}\right]\right)_{0}\cong\mathcal{I}_{F}\mathcal{RP}_{+}(F)\left[\frac{1}{2}\right]\rightarrow
\displaystyle\bigoplus_{p\in\Pi}\mathcal{P}(k(\upsilon_{p}))\left[\frac{1}{2}\right]\{p\}$
\end{center}

In \cite{ArticlehomologyQprint}, Hutchinson proved this is an isomorphism for $F=\mathbb{Q}$. Our main theorem (see Proposition \ref{Main_Proposition_Isomorphism}) is that it is an isomorphism for $F=\mathbb{Q}(\sqrt{-m})$ for $m\in \{2,3,,7,11\}$.
Furthermore we will prove a slightly weaker statement for the case $m=-1$. In all cases, we will use the local-global
principle for characters; i.e. we will prove that

\begin{center}
$\overline{S}_{\chi}:H_{0}\left(\mathcal{O}^{\times}_{F}, \mathcal{I}_{F}\mathcal{RP}_{+}(F)\left[\frac{1}{2}\right]\right)_{\chi}\rightarrow
\displaystyle\bigoplus_{p\in\Pi}\left(\mathcal{P}(k(\upsilon_{p}))\left[\frac{1}{2}\right]\{p\}\right)_{\chi}$
\end{center}

is an isomorphism for all $\chi\in \widehat{F^{\times}/(F^{\times})^{2}}$.

\smallskip

We will use the following observation repeatedly below:

\begin{lema}\label{Lemma_Action_of_square_classes_different_sign_characters}
Let $M$ be an $R_{F}$-module. Let $\chi\in\widehat{F^{\times}/(F^{\times})^{2}}$, let $\varepsilon\in\{-1, 1\}$. Suppose there exist $a\in F^{\times}$ such that $\langle a\rangle$ acts as $\varepsilon$ on $M$ and $\chi(\langle a\rangle)=-\varepsilon$. Then $M\left[\frac{1}{2}\right]_{\chi}=0$
\end{lema}

$\textbf{\emph{Proof}}.$ Let $m\in M$. In $M\left[\frac{1}{2}\right]_{\chi}$ we get

\begin{center}
$-\varepsilon m_{\chi}=\chi(\langle a\rangle) m_{\chi}=\langle a\rangle m_{\chi}=\varepsilon m_{\chi}$.
\end{center}

Then $2\varepsilon m_{\chi}=0$. Therefore $m_{\chi}=0$ for all $m_{\chi}\in M\left[\frac{1}{2}\right]_{\chi}$.

\begin{flushright}
    $\Box$
\end{flushright}

To start with, we consider the trivial character $\chi_{0}$. By definition $\mathcal{I}_{F}\mathcal{RP}_{+}(F)$ is generated by elements of the form
$\langle\langle x\rangle\rangle[y]$. Note that the square class of $x$ acts as multiplication by $-1$  on such an element (Since $\langle
x\rangle\langle\langle x\rangle\rangle=-\langle\langle x\rangle\rangle$). On the other hand, $\langle x\rangle$
acts as $1$ on any element of $\mathcal{I}_{F}\mathcal{RP}_{+}(F)_{\chi_{0}}$. By Lemma \ref{Lemma_Action_of_square_classes_different_sign_characters} we get $\left(\mathcal{I}_{F}\mathcal{RP}_{+}(F)\left[\frac{1}{2}\right]\right)_{\chi_{0}}=0$.


\smallskip

Similarly for $\mathcal{P}(k(\upsilon_{p}))\{p\}$, the element
$\langle p\rangle$ acts as $-1$ (see Remark \ref{squareclassuniformizer}) and
since $\chi_{0}(p)=1$, by Lemma \ref{Lemma_Action_of_square_classes_different_sign_characters} we deduce that
$\left(\mathcal{P}(k(\upsilon_{p}))\left[\frac{1}{2}\right]\{p\}\right)_{\chi_{0}}=0$.

\bigskip

Next, we consider characters $\chi$ satisfying $\chi(-1)=-1$. Let us recall that $\langle-1\rangle$ acts as multiplication by $1$ on
$\mathcal{I}_{F}\mathcal{RP}_{+}(F)\left[\frac{1}{2}\right]$ by Proposition
\ref{PropsitionRP+(F)andH3(SL2)} and also $\langle-1\rangle$ acts trivially
on $\mathcal{P}(k(\upsilon_{p}))$, since $-1$ is a unit for all $p\in\Pi$. By Lemma \ref{Lemma_Action_of_square_classes_different_sign_characters}, it follows that

\begin{center}
$\left(\mathcal{I}_{F}\mathcal{RP}_{+}(F)\left[\frac{1}{2}\right]\right)_{\chi}=0=
\displaystyle\bigoplus_{p\in\Pi}\left(\mathcal{P}(k(\upsilon_{p}))\left[\frac{1}{2}\right]\{p\}\right)_{\chi}$.
\end{center}

Hence we only need to consider $\chi\in
\widehat{F^{\times}/(F^{\times})^{2}}$ such that $\chi\not=\chi_{0}$ and
$\chi(-1)=1$.



\begin{lema}\label{LemmaTrivialActionOfSquareClass}
Let $F=\mathbb{Q}(\sqrt{-m})$ be an imaginary quadratic number field with $m>0$ and
$m\not=1$. Also let $\mathcal{O}_{F}$ be the ring of integers of $F$ then
$\langle u\rangle$ acts trivially on
$H_{3}\left(\mathrm{SL}_{2}(F),\mathbb{Z}\left[\frac{1}{2}\right]\right)$ for
all $u\in U=(\mathcal{O}_{F})^{\times}$.
\end{lema}

$\textbf{\emph{Proof}}.$ We have the following cases:

\begin{itemize}
\item[a)] If $m\not=-3$, then $U=(\mathcal{O}_{F})^{\times}=\{1,-1\}$.
    Since $\langle -1\rangle$ acts trivially on
    $H_{3}\left(\mathrm{SL}_{2}(F),\mathbb{Z}\left[\frac{1}{2}\right]\right)$,
    the result follows.
\item[b)] If $m=-3$ , then
    $U=(\mathcal{O}_{F})^{\times}=\{\pm1,\pm\omega,\pm\omega^{2}\}$ where
    $\omega^{3}=1$. It follows that $\omega=\omega^{4}$ and therefore
    $\langle\omega\rangle=\langle1\rangle$. Clearly we have that
    $\langle\omega^{2}\rangle=\langle1\rangle$ and therefore we get our
    result.

\end{itemize}

\begin{flushright}
    $\Box$
\end{flushright}





Now let $F=\mathbb{Q}(\sqrt{-m})$ be an imaginary quadratic number field with
$m\not=1$ and let $\chi\in \widehat{F^{\times}/(F^{\times})^{2}}$. If
$\chi(u)=-1$ for some $u\in (\mathcal{O}_{F})^{\times}$, then by Lemma \ref{Lemma_Action_of_square_classes_different_sign_characters} and Lemma \ref{LemmaTrivialActionOfSquareClass}, we get $\mathcal{RP}_{+}(F)\left[\frac{1}{2}\right]_{\chi}=0$. Furthermore, by Lemma \ref{Lemma_Action_of_square_classes_different_sign_characters} we also deduce $\left(\mathcal{P}(k(\upsilon_{p}))\left[\frac{1}{2}\right]\{p\}\right)_{\chi}=0$. 


\bigskip


If the class number of $F$ is $1$ then each discrete valuation $\upsilon$ is associated to a prime $p$ and

\begin{center}
$\frac{F^{\times}}{(F^{\times})^{2}(\mathcal{O}_{F})^{\times}}\cong\displaystyle\bigoplus_{p\in\Pi}p^{\mathbb{Z}/2}$.
\end{center}

(In fact, if $F$ has odd class number this latter isomorphism still holds, if for each prime ideal $\mathcal{P}$
we choose $\mathcal{P}^{m}=<p>$ for some odd $m$.)

\smallskip

It follows that

\begin{center}
$\widehat{\frac{F^{\times}}{(F^{\times})^{2}(\mathcal{O}_{F})^{\times}}}\cong\widehat{\displaystyle\bigoplus_{p\in\Pi}p^{\mathbb{Z}/2}}=\mathrm{Hom}_{\mathbb{Z}}\left(\displaystyle\bigoplus_{p\in\Pi}p^{\mathbb{Z}/2}, \mu_{2}\right)=\mathrm{Hom}_{\mathbf{Sets}}(\Pi,\mu_{2})$.
\end{center}

Thus
$\widehat{\frac{F^{\times}}{(F^{\times})^{2}(\mathcal{O}_{F})^{\times}}}$ is
naturally parameterized by the subsets of the set $\Pi$: If $S\in\Pi$ then
the corresponding character $\chi_{S}$ is defined by

\begin{center}
$\chi_{S}(p)=\left\{
               \begin{array}{ll}
                 -1, & \hbox{$p\in S$;} \\
                 1, & \hbox{$p\not\in S$.}
               \end{array}
             \right.
$
\end{center}

for all $p\in \Pi$ or, equivalently, for $x\in F^{\times}$, letting $\mathcal{L}(x):=\displaystyle\sum_{p\in S}\upsilon_{p}(x)$

\begin{center}
$\chi_{S}(x)=(-1)^{\mathcal{L}(x)}$
\end{center}

for all $x\in F^{\times}$. Conversely, the character $\chi$ is equal to $\chi_{S}$ where

\begin{center}
$S:=\mathrm{Supp}(\chi):=\{p\in \Pi \;|\; \chi(p)=-1\}$.
\end{center}

Thus, for a prime element $p$ in $\mathcal{O}_{F}$, $\chi_{p}$ is the unique
character satisfying $\mathrm{Supp}(\chi_{p})=\{p\}$.

\bigskip


We consider three kinds of characters:

\begin{itemize}
\item $\chi=\chi_{0}$ ($\mathrm{Supp}(\chi)=\emptyset$).
\item $\chi=\chi_{p}$ for some $p\in \mathcal{O}_{F}$; i.e.
    $|\mathrm{Supp}(\chi)|=1$
\item  $|\mathrm{Supp}(\chi)|\geq2$.
\end{itemize}


The following lemmas are immediate from the definition of the
$R_{F}$-structure on $\mathcal{P}(k(\upsilon_{p}))\{p\}$.

\begin{lema}\label{LemmaQuotientCharacterBlochgroup}Let $\chi\in\widehat{F^{\times}/(F^{\times})^{2}}$. Let $p$ be a prime element. Then

\begin{center}
$\left(\mathcal{P}(k(\upsilon_{p}))\left[\frac{1}{2}\right]\{p\}\right)_{\chi}=\left\{
                                                         \begin{array}{ll}
                                                           \mathcal{P}(k(\upsilon_{p}))\left[\frac{1}{2}\right]\{p\}, & \hbox{$\chi=\chi_{p}$;} \\
                                                           0, & \hbox{$otherwise$.}
                                                         \end{array}
                                                       \right.
$
\end{center}

\end{lema}

$\textbf{\emph{Proof}}.$ If $\chi\not=\chi_{p}$, then there exist $\langle
a\rangle\in F^{\times}/(F^{\times})^{2}$ such that $\chi(\langle
a\rangle)\not=\chi_{p}(\langle a\rangle)$; i.e. $\chi(\langle
a\rangle)=-\chi_{p}(\langle a\rangle)$. Therefore, the element $\langle
a\rangle$ acts in
$\left(\mathcal{P}(k(\upsilon_{p}))\left[\frac{1}{2}\right]\{p\}\right)_{\chi}$
as multiplication by $\chi(\langle a\rangle)$ and multiplication by
$\chi_{p}(\langle a\rangle)=-\chi(\langle a\rangle)$ . By Lemma \ref{Lemma_Action_of_square_classes_different_sign_characters} we get
$\left(\mathcal{P}(k(\upsilon_{p}))\left[\frac{1}{2}\right]\{p\}\right)_{\chi}=0$.

\smallskip

On other hand, when $\chi=\chi_{p}$, we have $\mathcal{P}(k(\upsilon_{p}))\{p\}=\mathcal{R}_{\chi_{p}}\otimes\mathcal{P}(k(\upsilon_{p}))= \mathcal{P}(k(\upsilon_{p}))_{\chi_{p}}$  by definition.

\begin{flushright}
    $\Box$
\end{flushright}

\begin{cor}\label{CorollaryQuotientCharacterBlochgroupDirectSum}
For $\chi\in\widehat{F^{\times}/(F^{\times})^{2}}$, we get

\begin{center}
$\left(\displaystyle\bigoplus_{p\in\Pi}\mathcal{P}(k(\upsilon_{p}))\left[\frac{1}{2}\right]\{p\}\right)_{\chi}=\left\{
                                                         \begin{array}{ll}
                                                           \mathcal{P}(k(\upsilon_{p}))\left[\frac{1}{2}\right]\{p\}, & \hbox{$\chi=\chi_{p}\; for\; some\; prime\; element \; p$;} \\
                                                           0, & \hbox{$otherwise$.}
                                                         \end{array}
                                                       \right.
$
\end{center}
\end{cor}

\begin{flushright}
    $\Box$
\end{flushright}







The following lemmas will allow us to prove that $\mathcal{RP}_{+}(F)\left[\frac{1}{2}\right]_{\chi}=0$ when $|\mathrm{Supp}(\chi)|\geq2$
(see Proposition \ref{Proposition_RP+(F)=0_Quotient_Characters})


\begin{lema}\label{Z[omega_m]contained_in_sum}
Let $F=\mathbb{Q}(\sqrt{m})$ be a quadratic number field, let $\alpha\not\sim\beta$ be prime elements of
$\mathcal{O}_{F}=\mathbb{Z}[\omega_{m}]$. Also let
$x\in\{\alpha, \alpha^{-1}\}$, $y\in\{\omega_{m}\alpha,
(\omega_{m}\alpha)^{-1}\}$, $z\in\{\beta,\beta^{-1}\}$ and
$w\in\{\omega_{m}\beta,(\omega_{m}\beta)^{-1}\}$. Then:

\begin{center}
$\mathcal{O}_{F}\subseteq
x\mathbb{Z}[x]+y\mathbb{Z}[y]+z\mathbb{Z}[z]+w\mathbb{Z}[w]$.
\end{center}

\end{lema}

$\textbf{\emph{Proof}}.$ Note by Lemma \ref{corollary_1_isintheset} and
Remark \ref{Remark_for_corollary_l_isintheset}, $\ell\mathbb{Z}[\ell]\subseteq
\frac{1}{\ell}\mathbb{Z}\left[\frac{1}{\ell}\right]$  for all non zero
$\ell\in \mathbb{Z}[\omega]$. So it is enough to prove the result for the
case when $x=\alpha$, $y=\omega_{m}\alpha$, $z=\beta$ and
$w=\omega_{m}\beta$. We have the following cases:

\begin{itemize}
\item[a)]  $\alpha=p\in \mathbb{Z}$ and $\beta=q\in \mathbb{Z}$; i.e.
    $\alpha$ and $\beta$ are primes in $\mathbb{Z}$ that inert in
    $\mathcal{O}_{F}$.

\smallskip

Note that $p\mathbb{Z}[p]=p\mathbb{Z}$. Thus

\begin{eqnarray}
\nonumber p\mathbb{Z}[p]+ q\mathbb{Z}[q] &=&p\mathbb{Z}+q\mathbb{Z}\\
\nonumber &=&\mathbb{Z}.
\end{eqnarray}

Also note that
$\omega_{m}p\mathbb{Z}[\omega_{m}p]=p^{2}N(\omega_{m})\mathbb{Z}+p\omega_{m}\mathbb{Z}$
from Lemma \ref{Corollary_lZ[l]_NumberField}. Therefore

\begin{eqnarray}
\nonumber \omega_{m}p\mathbb{Z}[\omega_{m}p]+ \omega_{m}q\mathbb{Z}[\omega_{m}q]&=&p^{2}N(\omega_{m})\mathbb{Z}+p\omega_{m}\mathbb{Z}+q^{2}N(\omega_{m})\mathbb{Z}+q\omega_{m}\mathbb{Z}\\
\nonumber &\supset& p\omega_{m}\mathbb{Z}+ q\omega_{m}\mathbb{Z}\\
\nonumber &=&\omega_{m}\mathbb{Z}
\end{eqnarray}

Hence

\begin{eqnarray}
\nonumber p\mathbb{Z}[p]+ q\mathbb{Z}[q]+\omega_{m}p\mathbb{Z}[\omega_{m}p]+ \omega_{m}q\mathbb{Z}[\omega_{m}q] &\supset& \mathbb{Z}+\omega_{m}\mathbb{Z}\\
\nonumber &=&\mathcal{O}_{F}.
\end{eqnarray}

\item[b)] $N(\alpha)$ is a prime number and $\beta=p\in \mathbb{Z}$; i.e.
    $\beta$ is a prime in $\mathbb{Z}$ that inert in
    $\mathcal{O}_{F}$

\smallskip

Let $\alpha=R(\alpha)+\omega_{m}I(\alpha)$. Note that

\begin{eqnarray}
\nonumber \alpha\mathbb{Z}[\alpha]&=&N(\alpha)\mathbb{Z}+\alpha\mathbb{Z}\\
\nonumber &=&\alpha(\overline{\alpha}\mathbb{Z}+\mathbb{Z})\\
\nonumber &=&\alpha(\mathbb{Z}+I(\overline{\alpha})\omega_{m}\mathbb{Z})\\
\nonumber &\supset&\alpha\mathbb{Z}.
\end{eqnarray}

Changing $\alpha$ by $\omega_{m}\alpha$, we get that

\begin{eqnarray}
\nonumber \omega_{m}\alpha\mathbb{Z}[\omega_{m}\alpha]&=&N(\omega_{m})N(\overline{\alpha})\mathbb{Z}+\alpha\omega_{m}\mathbb{Z}\\
\nonumber &=&\alpha(N(\omega_{m})\overline{\alpha}\mathbb{Z}+\omega_{m}\mathbb{Z})\\
\nonumber &=&\alpha(\mathrm{R}(\overline{\alpha})N(\omega_{m})\mathbb{Z}+I(\overline{\alpha})N(\omega_{m})\omega_{m}\mathbb{Z}+\omega_{m}\mathbb{Z})\\
\nonumber &\supset&\alpha\omega_{m}\mathbb{Z}. 
\end{eqnarray}

Therefore

\begin{eqnarray}
\nonumber \alpha\mathbb{Z}[\alpha]+\omega_{m}\alpha\mathbb{Z}[\omega_{m}\alpha]&\supset&\alpha\mathbb{Z}+\alpha\omega_{m}\mathbb{Z}\\
\nonumber &=&\alpha\mathcal{O}_{F}
\end{eqnarray}

Since $p\mathbb{Z}[p]=p\mathbb{Z}$. We get that

\begin{eqnarray}
\nonumber p\mathbb{Z}[p]+ \omega_{m}p\mathbb{Z}[\omega_{m}p]&=&p\mathbb{Z}+ p^{2}N(\omega_{m})\mathbb{Z}+\omega_{m}p\mathbb{Z}\\
\nonumber &\supset&p\mathbb{Z}+\omega_{m}p\mathbb{Z}\\
\nonumber &=&p(\mathbb{Z}+\omega_{m}\mathbb{Z})\\
\nonumber &=&p\mathcal{O}_{F}.
\end{eqnarray}

Therefore

\begin{eqnarray}
\nonumber p\mathbb{Z}[p]+ \omega_{m}p\mathbb{Z}[\omega_{m}p]+\alpha\mathbb{Z}[\alpha]+\omega_{m}\alpha\mathbb{Z}[\omega_{m}\alpha]&\supset&p\mathcal{O}_{F}+\alpha\mathcal{O}_{F}\\
\nonumber &=&\langle p, \alpha\rangle\\
\nonumber &=&\mathcal{O}_{F}.
\end{eqnarray}

\item[c)] $N(\alpha), N(\beta)$ are prime integers.

\smallskip

From case $b)$ we get that:

\begin{center}
$\alpha\mathbb{Z}[\alpha]+\omega_{m}\alpha\mathbb{Z}[\omega_{m}\alpha]\supset
\alpha\mathcal{O}_{F}$
\end{center}

and

\begin{center}
$\beta\mathbb{Z}[\beta]+\omega_{m}\beta\mathbb{Z}[\omega_{m}\beta]\supset
\beta\mathcal{O}_{F}$
\end{center}

Therefore

\begin{eqnarray}
\nonumber \alpha\mathbb{Z}[\alpha]+\omega_{m}\alpha\mathbb{Z}[\omega_{m}\alpha]+ \beta\mathbb{Z}[\beta]+\omega_{m}\beta\mathbb{Z}[\omega_{m}\beta]&\supset&\alpha\mathcal{O}_{F}+\beta\mathcal{O}_{F}\\
\nonumber &=&\langle\alpha,\beta\rangle\\
\nonumber &=&\mathcal{O}_{F}.
\end{eqnarray}

\end{itemize}

\begin{flushright}
    $\Box$
\end{flushright}

\begin{lema}\label{Z[omega_m]contained_in_sum for_X(omega)=-1}
Let $F$ be an imaginary quadratic number field, let $\alpha\not\sim\beta$ be prime elements of
$\mathcal{O}_{F}$. Also let
$x\in\{\alpha,\alpha^{-1}\}$, $y\in\{\beta,\beta^{-1}\}$ and
$z\in\{\omega_{m},\omega_{m}^{-1}\}$. Then:

\begin{center}
$\mathcal{O}_{F}\subseteq
x\mathbb{Z}[x]+y\mathbb{Z}[y]+z\mathbb{Z}[z]$.
\end{center}

\end{lema}

$\textbf{\emph{Proof}}.$ Similarly to Lemma \ref{Z[omega_m]contained_in_sum},
by Lemma \ref{corollary_1_isintheset} and Remark
\ref{Remark_for_corollary_l_isintheset}, $\ell\mathbb{Z}[\ell]\subseteq
\frac{1}{\ell}\mathbb{Z}\left[\frac{1}{\ell}\right]$  for all non zero
$\ell\in \mathbb{Z}[\omega]$. So it is enough to prove the result for the
case when $x=\alpha$, $y=\beta$ and $z=\omega_{m}$. We have the following
cases:

\begin{itemize}
\item[a)]  $\alpha=p\in \mathbb{Z}$ and $\beta=q\in \mathbb{Z}$; i.e.
    $\alpha$ and $\beta$ are primes in $\mathbb{Z}$ that inert in
    $\mathcal{O}_{F}$.

\smallskip

Note that $p\mathbb{Z}[p]=p\mathbb{Z}$. It follows that
$p\mathbb{Z}[p]+q\mathbb{Z}[q]=p\mathbb{Z}+q\mathbb{Z}=\mathbb{Z}$. Hence

\begin{eqnarray}
\nonumber p\mathbb{Z}[p]+q\mathbb{Z}[q]+ \omega_{m}\mathbb{Z}[\omega_{m}]&=&\mathbb{Z}+N(\omega_{m})\mathbb{Z}+\omega_{m}\mathbb{Z}\\
\nonumber &\supset&\mathbb{Z}+ \omega_{m}\mathbb{Z}\\
\nonumber &=&\mathcal{O}_{F}.
\end{eqnarray}

\item[b)] $N(\alpha)$ is a prime number and $\beta=p\in \mathbb{Z}$; i.e.
    $\beta$ is a prime in $\mathbb{Z}$ that inert in
    $\mathcal{O}_{F}$

\smallskip

Let $\alpha=R(\alpha)+\omega_{m}I(\alpha)$. Since $\alpha\not\sim
\beta=p$, it follows that $N(\alpha)\not=p$ and thus
$N(\alpha)\mathbb{Z}+p\mathbb{Z}=\mathbb{Z}$. Therefore

\begin{eqnarray}
\nonumber \alpha\mathbb{Z}[\alpha]+\beta\mathbb{Z}[\beta]+\omega_{m}\mathbb{Z}[\omega_{m}]&=&N(\alpha)\mathbb{Z}+\alpha\mathbb{Z}+p\mathbb{Z}+N(\omega_{m})\mathbb{Z}+\omega_{m}\mathbb{Z}\\
\nonumber &=&\left(N(\alpha)\mathbb{Z}+p\mathbb{Z}\right)+N(\omega_{m})\mathbb{Z}+\alpha\mathbb{Z}+\omega_{m}\mathbb{Z}\\
\nonumber &=&\mathbb{Z}+N(\omega_{m})\mathbb{Z}+\alpha\mathbb{Z}+\omega_{m}\mathbb{Z}\\
\nonumber &\supset&\mathbb{Z}+\omega_{m}\mathbb{Z}\\
\nonumber &=&\mathcal{O}_{F}.
\end{eqnarray}

\item[c)] $N(\alpha), N(\beta)$ are prime integers.

\smallskip

We will consider the following cases:

\begin{itemize}

\item[1)] $N(\alpha)\not=N(\beta)$.

Then

\begin{eqnarray}
\nonumber \alpha\mathbb{Z}[\alpha]+\beta\mathbb{Z}[\beta]&=& N(\alpha)\mathbb{Z}+\alpha
\mathbb{Z}+N(\beta)\mathbb{Z}+\beta\mathbb{Z}\\
\nonumber &=&\mathbb{Z}+ \alpha
\mathbb{Z}+\beta\mathbb{Z}\\
\nonumber &\supset&\mathbb{Z}
\end{eqnarray}

Therefore

\begin{eqnarray}
\nonumber \alpha\mathbb{Z}[\alpha]+\beta\mathbb{Z}[\beta]+ \omega_{m}\mathbb{Z}[\omega_{m}]&\supset&\mathbb{Z}+ N(\omega_{m})\mathbb{Z}+\omega_{m}\mathbb{Z}\\
\nonumber &\supset&\mathbb{Z}+\omega_{m}\mathbb{Z}\\
\nonumber &=&\mathcal{O}_{F}.
\end{eqnarray}

\item[2)] $N(\alpha)=N(\beta)=p$.

Then for $\alpha=R(\alpha)+I(\alpha)\omega_{m}$, we have
$\beta=\overline{\alpha}$. We consider the following subcases:

\begin{itemize}
\item[a)] $m\;\equiv\;1\;(\mathrm{mod}\;4)$.

Then $\beta=(R(\alpha)+I(\alpha))-I(\alpha)\omega_{m}$. Thus

\begin{eqnarray}
\nonumber \alpha\mathbb{Z}[\alpha]+\beta\mathbb{Z}[\beta]&=&N(\alpha)\mathbb{Z}+\alpha\mathbb{Z}+N(\beta)\mathbb{Z}+\beta\mathbb{Z}\\
\nonumber &=& N(\alpha)\mathbb{Z}+\alpha\mathbb{Z}+\beta\mathbb{Z}  \;\;\;\mbox{(Since $N(\alpha)=N(\beta)$)}\\
\nonumber &=&N(\alpha)\mathbb{Z}+ R(\alpha)\mathbb{Z}+ I(\alpha)\mathbb{Z}+ I(\alpha)\omega_{m}\mathbb{Z}\\
\nonumber &\supset&\mathbb{Z} + I(\alpha)\omega_{m}\mathbb{Z}\;\;\; \mbox{(Since $p\not|R(\alpha)$ or $p\not|I(\alpha)$ )}\\
\nonumber &\supset&\mathbb{Z}.
\end{eqnarray}

Therefore

\begin{eqnarray}
\nonumber \alpha\mathbb{Z}[\alpha]+\beta\mathbb{Z}[\beta]+ \omega_{m}\mathbb{Z}[\omega_{m}]&\supset&\mathbb{Z}+ N(\omega_{m})\mathbb{Z}+\omega_{m}\mathbb{Z}\\
\nonumber &\supset&\mathbb{Z}+\omega_{m}\mathbb{Z}\\
\nonumber &=&\mathcal{O}_{F}.
\end{eqnarray}

\item[b)] $m\;\equiv\;2,3\;(\mathrm{mod}\;4)$.

Then $\beta=R(\alpha)-I(\alpha)\omega_{m}$. Hence

\begin{eqnarray}
\nonumber \alpha\mathbb{Z}[\alpha]+\beta\mathbb{Z}[\beta]&=&N(\alpha)\mathbb{Z}+\alpha\mathbb{Z}+N(\beta)\mathbb{Z}+\beta\mathbb{Z}\\
\nonumber &=& N(\alpha)\mathbb{Z}+\alpha\mathbb{Z}+\beta\mathbb{Z}  \;\;\;\mbox{(Since $N(\alpha)=N(\beta)$)}\\
\nonumber &=&N(\alpha)\mathbb{Z}+ R(\alpha)\mathbb{Z}+ I(\alpha)\omega_{m}\mathbb{Z}.
\end{eqnarray}

Note that either $N(\alpha)$ and $R(\alpha)$ are coprimes or $N(\alpha)$ and $N(\omega_{m})$ are coprimes otherwise
by Lemma \ref{R(alpha)_and__N(alpha)_are_not_coprimes}, we get $\alpha\sim\beta$. Hence it follows:



\begin{eqnarray}
\nonumber \alpha\mathbb{Z}[\alpha]+\beta\mathbb{Z}[\beta]+ \omega_{m}\mathbb{Z}[\omega_{m}]&=&N(\alpha)\mathbb{Z}+ R(\alpha)\mathbb{Z}+ I(\alpha)\omega_{m}\mathbb{Z} + N(\omega_{m})\mathbb{Z}+\omega_{m}\mathbb{Z}\\
\nonumber &\supset&\mathbb{Z}+\omega_{m}\mathbb{Z}\\
\nonumber &=&\mathcal{O}_{F}.
\end{eqnarray}

\end{itemize}

\end{itemize}

\end{itemize}

\begin{flushright}
    $\Box$
\end{flushright}

\begin{prop}\label{Propositionsumalgebraiccintegergaussian}
Let $F$ be an imaginary quadratic number field of class number $1$ 
and let $\chi\in\widehat{\frac{F^{\times}}{(F^{\times})^{2}(\mathcal{O}_{F})^{\times}}}$.
Suppose that $|\mathrm{Supp}(\chi)|\geq2$. Then

\begin{center}
$[a]_{\chi}=[a+t]_{\chi}$
\end{center}

in $\mathcal{RP}_{+}(F)\left[\frac{1}{2}\right]_{\chi}$
for all $a\in F$  and $ t\in\mathcal{O}_{F}$.

\end{prop}

$\textbf{\emph{Proof}}.$ Let $p,q\in \mathrm{Supp}(\chi)$. We consider the
following cases:

\begin{itemize}

\item[a)] $\chi(\omega_{m})=1$.

\smallskip

Then $\chi(\omega_{m}p)=\chi(\omega_{m}q)=-1$. Therefore, depending of the values of $\chi(1-p)$, $\chi(1-q)$, $\chi(1-\omega_{m}p)$ and
$\chi(1-\omega_{m}q)$, by Remark \ref{Remark[a]=[la]}, there is a choice
of $\epsilon_{p},\epsilon_{q},\epsilon'_{p},\epsilon'_{q}\in \{\pm 1\}$
such that $\ell=p^{\epsilon_{p}}, q^{\epsilon_{q}},
(\omega_{m}p)^{\epsilon'_{p}},(\omega_{m}q)^{\epsilon'_{q}}$ satisfy
$\chi(\ell)=-1$ and $\chi(1-\ell)=1$. Therefore the result follows from
the Proposition \ref{Propositionelement_lZ[l]} and the Lemma
\ref{Z[omega_m]contained_in_sum}.

\item[b)] $\chi(\omega_{m})=-1$.

\smallskip

Then depending of the values of $\chi(1-p)$, $\chi(1-q)$ and
$\chi(1-\omega_{m})$, by the Remark \ref{Remark[a]=[la]} there is a
choice of $\epsilon_{p},\epsilon_{q},\epsilon_{\omega_{m}}\in \{\pm 1\}$ such that
$\ell=p^{\epsilon_{p}}, q^{\epsilon_{q}},
(\omega_{m})^{\epsilon_{\omega_{m}}}$ satisfy $\chi(\ell)=-1$ and
$\chi(1-\ell)=1$. Therefore from the Proposition
\ref{Propositionelement_lZ[l]} and the Lemma
\ref{Z[omega_m]contained_in_sum for_X(omega)=-1}, we get our result.

\end{itemize}

\begin{flushright}
    $\Box$
\end{flushright}

\begin{prop}\label{Proposition_RP+(F)=0_Quotient_Characters}
Let $F=\mathbb{Q}(\sqrt{-m})$ be an imaginary quadratic number field with $m=1,2,3,7,11$ and
$\chi\in\widehat{\frac{F^{\times}}{(F^{\times})^{2}(\mathcal{O}_{F})^{\times}}}$.
Suppose that $|\mathrm{Supp}(\chi)|\geq2$. Then
$\mathcal{RP}_{+}(F)\left[\frac{1}{2}\right]_{\chi}=0$.
\end{prop}

$\textbf{\emph{Proof}}.$ First, let us recall, for these fields, the ring of algebraic integer $\mathcal{O}_{F}$ is a Euclidean domain with respect the norm of $F$.

\smallskip

From the Proposition \ref{Propositionsumalgebraiccintegergaussian}, we get
$[x]_{\chi}=[x+a]_{\chi}$  for all $x\in F$ and $a\in \mathcal{O}_{F}$. Then
it follows that $[a]_{\chi}=[1]_{\chi}=0$ for all $a\in \mathcal{O}_{F}$ and
therefore $\left[\frac{1}{a}\right]_{\chi}=0$ for all $a\in \mathcal{O}_{F}$.
\bigskip

 Now for $x\in F$ we have that $x=\frac{a}{b}$ with $a,b \in \mathcal{O}_{F}$, we will prove that $[x]_{\chi}=0$ by induction on $\mathrm{min}(N(a),N(b))=m$. If $m=1$, then $a\in \mathcal{O}^{\times}_{F}$ or $b\in \mathcal{O}^{\times}_{F}$. Thus either $x\in \mathcal{O}_{F}$ or $x^{-1}\in \mathcal{O}_{F}$. Therefore, it follows from the previous paragraph that $[x]_{\chi}=0$.

\bigskip

Since $[x]_{\chi}=-\left[\frac{1}{x}\right]_{\chi}$, we can suppose that
$N(b)<N(a)$. Now we suppose the results holds if norm $m=N(b)\leq n$. We will
prove it for $\mathrm{min}(N(a),N(b))=N(b)=n+1$. Note that there exist
$q,r\in \mathcal{O}_{F}$ such that $a=qb+r$ and $N(r)<N(b)$. Taking $t=q$
from Proposition \ref{Propositionsumalgebraiccintegergaussian}, we obtain that
$[x]_{\chi}=\left[\frac{r}{b}\right]_{\chi}$. Note that
$\mathrm{min}(N(r),N(b))=N(r)\leq n$. Hence by induction hypothesis we have
that $[x]_{\chi}=0$.

\begin{flushright}
    $\Box$
\end{flushright}

\begin{lema}
If $f:M\rightarrow N$ is a map of $R_{F}$-modules and if $f_{\chi}:M_{\chi}\rightarrow N_{\chi}$ is an isomorphism for all $\chi\in\widehat{\frac{F^{\times}}{(F^{\times})^{2}(\mathcal{O}_{F})^{\times}}}$, then the induced map

\begin{center}
$H_{0}\left(\mathcal{O}^{\times}_{F}, M\left[\frac{1}{2}\right]\right)\rightarrow H_{0}\left(\mathcal{O}^{\times}_{F}, N\left[\frac{1}{2}\right]\right)$
\end{center}

is an isomorphism.
\end{lema}

$\textbf{\emph{Proof}}.$ Let $\chi\in\widehat{\frac{F^{\times}}{(F^{\times})^{2}}}$. If $\chi\in \widehat{\frac{F^{\times}}{(F^{\times})^{2}(\mathcal{O}_{F})^{\times}}}$, then $H_{0}\left(\mathcal{O}^{\times}_{F}, M\left[\frac{1}{2}\right]\right)_{\chi}=M\left[\frac{1}{2}\right]_{\chi}\cong N\left[\frac{1}{2}\right]_{\chi} =H_{0}\left(\mathcal{O}^{\times}_{F}, N\left[\frac{1}{2}\right]\right)_{\chi}$.

\bigskip

Otherwise, $\chi(u)=-1$ for some $u\in \mathcal{O}^{\times}$. But the element $\langle u\rangle$ acts trivially on $H_{0}\left(\mathcal{O}^{\times}_{F}, M\left[\frac{1}{2}\right]\right)$ and $H_{0}\left(\mathcal{O}^{\times}_{F}, N\left[\frac{1}{2}\right]\right)$. By Lemma \ref{Lemma_Action_of_square_classes_different_sign_characters}

\begin{center}
$H_{0}\left(\mathcal{O}^{\times}_{F}, M\left[\frac{1}{2}\right]\right)_{\chi}=0=H_{0}\left(\mathcal{O}^{\times}_{F}, N\left[\frac{1}{2}\right]\right)_{\chi}$.
\end{center}

The result follows by the local-global principle for characters (Proposition \ref{Proposition_Local_global_principle_character}).

\begin{flushright}
    $\Box$
\end{flushright}


\begin{prop}\label{Main_Proposition_Isomorphism}
Let $F=\mathbb{Q}(\sqrt{-m})$ be an imaginary quadratic number field with $m=1, 2, 3, 7, 11$ , then the map


\begin{center}
$\overline{S}:H_{0}\left(\mathcal{O}_{F}^{\times}, H_{3}\left(\mathrm{SL}_{2}(F),\mathbb{Z}\left[\frac{1}{2}\right]\right)_{0}\right)\cong H_{0}\left(\mathcal{O}_{F}^{\times},\mathcal{I}_{F}\mathcal{RP}_{+}(F)\left[\frac{1}{2}\right]\right)\rightarrow
\displaystyle\bigoplus_{p\in\Pi}\mathcal{P}(k(\upsilon_{p}))\left[\frac{1}{2}\right]\{p\}$
\end{center}

is an isomorphism of $R_{F}$-modules.

\end{prop}

$\textbf{\emph{Proof}}.$ By the local-global principle for characters
(Proposition \ref{Proposition_Local_global_principle_character}) and by
Corollary \ref{CorollaryQuotientCharacterBlochgroupDirectSum}, we must show
that

\begin{center}
$\left(\mathcal{I}_{F}\mathcal{RP}_{+}(F)\left[\frac{1}{2}\right]\right)_{\chi}=\left\{
                                                                                      \begin{array}{ll}
                                                                                        \mathcal{P}(k(\upsilon_{p}))\left[\frac{1}{2}\right]\{p\}, & \hbox{$if \; \chi=\chi_{p}$;} \\
                                                                                        0, & \hbox{$if\; \chi\not=\chi_{p}\; for\; any\; p$.}
                                                                                      \end{array}
                                                                                    \right.
$
\end{center}

If $\chi=\chi_{p}$, this is from Proposition \ref{NaturalIsomorphismRP+(F)}.
If $\chi=\chi_{0}$, or $\chi(u)=-1$ for some $u\in(\mathcal{O}_{F})^{\times}$
this follows from above. Otherwise
$\chi\in\widehat{\frac{F^{\times}}{(F^{\times})^{2}(\mathcal{O}_{F})^{\times}}}$
and  $\mathrm{Supp}(\chi)$ contain at least two non-associated distinct prime
elements of $\mathcal{O}_{F}$ and hence from Proposition
\ref{Proposition_RP+(F)=0_Quotient_Characters} we get our result.

\begin{flushright}
    $\Box$
\end{flushright}

From the Proposition \ref{Main_Proposition_Isomorphism} and the definition of
$H_{3}\left(\mathrm{SL}_{2}(F),\mathbb{Z}\left[\frac{1}{2}\right]\right)_{0}$,
we get:

\begin{cor}\label{Corollary_Short_exact_sequence_quadratic_extension}
Let $F=\mathbb{Q}(\sqrt{-m})$ be an imaginary quadratic number field with $m=1,2,3,7,11$ , then we have a short exact sequence of $R_{F}$-modules

\begin{center}
$\xymatrix{0\ar[r]&\displaystyle\bigoplus_{p\in\Pi}\mathcal{P}(k(\upsilon_{p}))\left[\frac{1}{2}\right]\{p\}\ar[r]&H_{0}\left(\mathcal{O}_{F}^{\times},H_{3}\left(\mathrm{SL}_{2}(F),\mathbb{Z}\left[\frac{1}{2}\right]\right)\right)\ar[r]&K_{3}^{\mathrm{ind}}(F)\left[\frac{1}{2}\right]\ar[r]&0}$.
\end{center}

\end{cor}


\begin{flushright}
    $\Box$
\end{flushright}

\begin{cor}
Let $F=\mathbb{Q}(\sqrt{-m})$ be an imaginary quadratic number field with $m=2,3,7,11$ ,then we have a short exact sequence of $R_{F}$-modules

\begin{center}
$\xymatrix{0\ar[r]&\displaystyle\bigoplus_{p\in\Pi}\mathcal{P}(k(\upsilon_{p}))\left[\frac{1}{2}\right]\{p\}\ar[r]&H_{3}\left(\mathrm{SL}_{2}(F),\mathbb{Z}\left[\frac{1}{2}\right]\right)\ar[r]&K_{3}^{\mathrm{ind}}(F)\left[\frac{1}{2}\right]\ar[r]&0}$.
\end{center}

\end{cor}

$\textbf{\emph{Proof}}.$ By Lemma \ref{LemmaTrivialActionOfSquareClass}, we have that $\mathcal{O}^{\times}_{F}$ acts trivially on $H_{3}\left(\mathrm{SL}_{2}(F),\mathbb{Z}\left[\frac{1}{2}\right]\right)$.

\begin{flushright}
    $\Box$
\end{flushright}

\end{document}